\documentclass[a4paper,11pt]{article}
\usepackage[top=2.5cm,bottom=2.5cm,left=2.2cm,right=2.2cm]{geometry}
\usepackage{amsfonts}
\usepackage{mathrsfs,amscd,amssymb,amsthm,amsmath,bm,graphicx,psfrag,subfigure,url,mathtools}
\usepackage{pict2e}
\usepackage{stfloats}
\usepackage{psfrag,amsmath}
\usepackage[T1]{fontenc}
\usepackage[utf8]{inputenc}
\usepackage{lmodern}
\usepackage{amsmath,amssymb,amsthm}
\usepackage{microtype}
\usepackage{enumitem}
\usepackage{booktabs}
\usepackage{tikz}
\usepackage{fancyhdr}
\usepackage{xcolor}
\usepackage{authblk}

\definecolor{linkblue}{HTML}{184A84}
\usepackage[hidelinks]{hyperref}
\hypersetup{pdftitle={Counting cliques in optimal 1-planar graphs},pdfauthor={}}
\theoremstyle{plain}
\newtheorem{theorem}{Theorem}[section]
\newtheorem{lemma}[theorem]{Lemma}
\newtheorem{proposition}[theorem]{Proposition}
\newtheorem{corollary}[theorem]{Corollary}
\theoremstyle{definition}
\newtheorem{question}{Question}

\newtheorem{construction}{Construction}
\newtheorem{conjecture}{Conjecture}
\theoremstyle{remark}

\DeclareMathOperator{\Rad}{Rad}

\setlist[enumerate]{itemsep=3pt,topsep=4pt,font=\normalfont}
\setlist[itemize]{itemsep=3pt,topsep=4pt}
\hypersetup{
    colorlinks=true,
    linkcolor=blue,
    citecolor=red,
    urlcolor=magenta,
}
\def\b{\color{blue}}

\newcommand{\keywords}[1]{%
  \par\vspace{6pt}\noindent\textbf{Keywords: }#1\par
}

\newcommand{\MSC}[2][2020]{%
  \par\vspace{3pt}\noindent\textbf{MSC(#1): }#2\par
}
\begin{document}
\baselineskip=0.21in
\title{\LARGE\bfseries Counting cliques in optimal $1$-planar graphs}
\author[1,2]{Shuchao Li}
\author[1]{Liwen Lian\thanks{Corresponding author: Liwen Lian. e-mail: lianliwen@mails.ccnu.edu.cn}}
\affil[1]{School of Mathematics and Statistics, and Hubei Key Lab--Math. Sci.,\linebreak Central China Normal University, Wuhan 430079, China}
\affil[2]{Key Laboratory of Nonlinear Analysis \& Applications (Ministry of Education),\linebreak Central China Normal University, Wuhan 430079, China}

\date{\today}
\maketitle
\begin{abstract}
Extremal problems concerning the maximum number of cliques in $n$-vertex graphs from a given graph class have been widely investigated. In this contribution, we present several results on counting cliques that extend beyond the existing framework. Firstly, we address this problem for the class of optimal 1-planar graphs where we determine the maximum numbers of copies of $K_3$, $K_4$, and $K_5$ in an $n$-vertex optimal 1-planar graphs and characterize all extremal graphs. This resolves Question {\b 8.1} of Gollin-Hendrey-Methuku-Tompkins-Zhang. 
Secondly, we prove that every 3-connected plane quadrangulation on $n\geq14$ vertices has at most $2n-14$ 3-separators and characterize the corresponding extremal graphs. This confirms Conjecture {\b 8.2} posed by Gollin-Hendrey-Methuku-Tompkins-Zhang. Lastly, we construct, for every integer $r\geq1$, a simple 2-planar graph with $5r+3$ vertices and $44r+1$ triangles, thus disproving Conjecture {\b 8.3} of theirs. Our proofs use radial graphs, tree decomposition and SPQR trees.
\end{abstract}
\MSC{05C10, 05C35, 05C40}

\section{Introduction}
\subsection{Background}
We consider finite undirected graphs, possibly with parallel edges
but without loops, unless stated otherwise.
For a graph $G$, let $V(G)$ and $E(G)$ denote its vertex set and
edge set, respectively.
Let $K_t$ denote the complete graph with $t$ vertices, and let
$C_n$ denote the cycle of length $n$.
For $X\subseteq V(G)$, let $G[X]$ denote the subgraph of $G$
induced by $X$.
For a simple graph $G$, a set $X\subseteq V(G)$ is a
\emph{clique} if $G[X]\cong K_t$ for some $t$.
For graphs $G$ and $H$, let $N(G,H)$ denote the number of
subgraphs of $G$ isomorphic to $H$.

For $X\subseteq V(G)$, let $G-X$ denote the graph obtained by
deleting the vertices in $X$ and their incident edges.
For a positive integer $k$, a graph $G$ is \emph{$k$-connected}
if $|V(G)|\ge k+1$ and $G-X$ is connected for every
$X\subseteq V(G)$ with $|X|<k$.
A set $X\subseteq V(G)$ is a \emph{$3$-separator} if $|X|=3$
and $G-X$ is disconnected.

A \emph{planar embedding} of a graph is a drawing in the plane
whose edges meet only at common ends.
A graph with a fixed planar embedding is a \emph{plane graph}.
The \emph{faces} of a plane graph are the connected components
of the complement of its drawing in the plane.
The unique unbounded face is the \emph{outer face}.
For a plane graph $G$, let $F(G)$ denote its set of faces,
and let $\partial f$ denote the boundary of a face $f\in F(G)$.
A \emph{plane quadrangulation} is a connected plane graph whose
face boundaries are cycles of length four.

A graph is \emph{$1$-planar} if it has a drawing in the plane
with each edge crossed at most once.
Ringel~\cite{ref01} introduced these graphs in 1965 in
connection with the simultaneous colouring of vertices and faces
of plane graphs.
He proved that every $1$-planar graph is $7$-colourable and
conjectured that six colours suffice.
Borodin~\cite{Borodin1984} proved this conjecture in 1984.
Surveys of structural, geometric and algorithmic results on
$1$-planar graphs can be found
in~\cite{ref04, ref05}.

Bodendiek, Schumacher and Wagner~\cite{ref02} proved
that every simple $1$-planar graph with $n\ge3$ vertices has
at most $4n-8$ edges.
A simple $1$-planar graph with $n\ge3$ vertices and exactly
$4n-8$ edges is called \emph{optimal $1$-planar}.
For a $3$-connected plane quadrangulation $Q$, let $\widehat Q$
denote the graph obtained by adding both diagonals in every face,
with the two diagonals crossing inside that face.
Bodendiek, Schumacher and Wagner also characterized optimal
$1$-planar graphs as precisely the graphs obtained in this way.
They proved that an optimal $1$-planar graph with $n$ vertices
exists if and only if $n=8$ or $n\ge10$
\cite{ref02,Bodendiek1984}.

The quadrangulation characterization implies that every vertex
of an optimal $1$-planar graph has even degree at least six.
Moreover, every optimal $1$-planar graph has at least eight
vertices of degree six~\cite{Brandenburg2018}.
Fujisawa, Segawa and Suzuki~\cite{FujisawaSegawaSuzuki2018}
proved that its vertex-connectivity is either four or six.
Noguchi and Suzuki~\cite{NoguchiSuzuki2015} showed that every
optimal $1$-planar graph contains a $4$-connected
triangulation as a spanning subgraph, and hence is Hamiltonian.
The $4$-connectivity, together with \cite[Lemma~6.1]{ref07} of
Gollin et al., also implies that optimal
$1$-planar graphs are $K_6$-free. Other work on optimal $1$-planar graphs concerns recognition algorithms, matching extension and embedding.
For details, see~\cite{Brandenburg2018,FujisawaSegawaSuzuki2018,Suzuki2010}.

Determining the maximum number of copies of a complete subgraph
in a given graph class is a basic problem in extremal graph theory.
For simple planar graphs with $n\ge3$ vertices, Euler's formula
gives the edge bound $3n-6$~\cite{ref03}.
Hakimi and Schmeichel~\cite{ref16} proved that
the maximum number of triangles is $3n-8$.
For $n\ge4$, Alon and Caro~\cite{ref17} and
Wood~\cite{ref11}, independently, proved that the maximum
number of copies of $K_4$ is $n-3$.
The graphs attaining either bound are precisely the
\emph{Apollonian networks}: these are obtained from
$K_3$ by repeatedly adding a vertex inside a triangular face
and joining it to the three boundary vertices.
Clique-counting problems have also been studied for graphs
embedded in a fixed surface~\cite{ref27} and for graphs
excluding a fixed minor~\cite{ref31,ref33}.

Huang and Zhang~\cite{HuangZhang2026} obtained
sharp upper bounds for the numbers of triangles and copies
of $K_4$ in $7$-connected $1$-planar graphs. Gollin et al.~\cite{ref07} determined the maximum number
of copies of $K_t$ in simple $1$-planar graphs for each
$t\in\{3,4,5,6\}$ and characterized all extremal graphs.
Their proofs combine planar clique bounds with edge deletions
and induction along separating triangles. For $n\ge9$ and $t\in\{3,4,5,6\}$, every simple $1$-planar
graph with $n$ vertices that maximizes $N(G,K_t)$ has fewer
than $4n-8$ edges.
Gollin et al.~\cite{ref07} therefore asked the following
question for optimal $1$-planar graphs.

\begin{question}[\cite{ref07}, Question~8.1]
\label{ques}
Let $n\ge10$ and $t\in\{3,4,5\}$ be integers.
What is the maximum number of copies of $K_t$ in an optimal
$1$-planar graph with $n$ vertices?
\end{question}

For a $3$-connected plane quadrangulation $Q$, let $s_3(Q)$
denote the number of its $3$-separators.
Gollin et al.~\cite{ref07} observed that, for $k=3$,
Question~\ref{ques} is equivalent to maximizing $s_3(Q)$ over
$3$-connected plane quadrangulations with $n$ vertices.
By Corollary~\ref{cor:triangles},
$N(\widehat Q,K_3)=4n-8+s_3(Q)$, where $n=|V(Q)|$.
Gollin et al.~\cite{ref07} constructed the following family
$\mathcal G^*$.

\begin{construction}
\label{con}
Let $H_0$ be the $4$-cycle $v_0v_1v_2v_3v_0$, and let
$\varphi_0$ be a planar embedding of $H_0$ with faces
$F_1^0$ and $F_2^0$.
For each $i\ge1$, obtain $H_i$ from $H_{i-1}$ by adding a
new vertex $v_{i+3}$ adjacent to $v_{i+2}$ and to the vertex
opposite $v_{i+2}$ on the boundary of $F_2^{i-1}$.
Extend $\varphi_{i-1}$ to a planar embedding $\varphi_i$
of $H_i$ by drawing the new vertex and edges inside
$F_2^{i-1}$.
Let $F_1^i=F_1^{i-1}$, and choose either of the two new
faces as $F_2^i$.
We call $F_2^i$ the \emph{active face}.

For each $i\ge1$, take two disjoint copies of the cube $Q_3$.
Identify a facial $4$-cycle in the first copy with the boundary
cycle of $F_1^i$, and a facial $4$-cycle in the second copy
with the boundary cycle of $F_2^i$, drawing the remaining
vertices and edges of each cube in the respective face.
Let $G_i$ denote the resulting graph.
The family $\mathcal G^*$ consists of all graphs $G_i$
obtained for $i\ge1$ and all choices of the active faces.
By construction, $|V(G_i)|=i+12$. Moreover, Every graph in $\mathcal G^*$ is a $3$-connected plane
quadrangulation~\cite{ref07}.
\end{construction}

Gollin et al.~\cite{ref07} conjectured that this family
gives precisely the extremal quadrangulations for $n\ge14$.

\begin{conjecture}[\cite{ref07}, Conjecture~8.2]
\label{conj}
For every integer $n\ge14$, the $3$-connected planar
quadrangulations with $n$ vertices and the maximum number
of $3$-separators are precisely the graphs in $\mathcal G^*$
with $n$ vertices.
\end{conjecture}

A graph is \emph{$2$-planar} if it has a drawing in the plane
with each edge crossed at most twice.
The complete graph $K_t$ is $2$-planar if and only if
$t\le7$~\cite{Angelini2020}.
For each $t\in\{3,4,5,6,7\}$, the exact maximum number
of copies of $K_t$ in a simple $2$-planar graph
with $n$ vertices remains unknown.
Gollin et al.~\cite{ref07} proposed the following conjecture on the largest number of triangles in simple $2$-planar graphs.
\begin{conjecture}[\cite{ref07}, Conjecture~8.3]
\label{conj2}
Every simple $2$-planar graph $G$ with $n\ge7$ vertices
satisfies
\(
N(G,K_3)\le \frac{17n-49}{2}.
\)
\end{conjecture}

\subsection{Main results}

In this paper, we mainly focus on clique counting problems for optimal 1-planar graphs. We answer Question~\ref{ques},
prove Conjecture~\ref{conj}, and disprove
Conjecture~\ref{conj2}. 

Our first result establishes the sharp upper bound on the number
of $3$-separators and characterizes all extremal quadrangulations, which confirms Conjecture~\ref{conj}.
\begin{theorem}
\label{thm}
Let $Q$ be a $3$-connected plane quadrangulation with
$n\ge14$ vertices.
Then
\[
s_3(Q)\le 2n-14,
\]
with equality if and only if $Q\in\mathcal G^{*}$,
where $\mathcal G^{*}$ is defined in
Construction~\ref{con}.
\end{theorem}

For a connected plane graph $H$, let $\operatorname{Rad}(H)$
denote the bipartite plane graph with partite sets $V(H)$ and
$F(H)$, where $v\in V(H)$ is adjacent to $f\in F(H)$ if and
only if $v$ lies on the boundary of $f$.
For $n\ge3$, let $W_n$ denote the wheel with rim $C_n$.
The triangular prism is $C_3\square K_2$, where $\square$
denotes the Cartesian product.
Our next theorem answers Question~\ref{ques}.
\begin{theorem}\label{thm:cliques}
For each integer $n\ge10$, the maximum values of $N(G,K_k)$
over all optimal $1$-planar graphs $G$ with $n$ vertices,
for $k\in\{3,4,5\}$, are given in the following table.

\begin{center}
\renewcommand{\arraystretch}{1.12}
\begin{tabular}{cccc}
\toprule
$n$ & $K_3$ & $K_4$ & $K_5$ \\
\midrule
$10$ & $40$ & $16$ & $0$ \\
$11$ & $45$ & $19$ & $0$ \\
$12$ & $50$ & $20$ & $0$ \\
$13$ & $57$ & $25$ & $0$ \\
$n\ge14$ & $6n-22$ & $4n-29$ & $n-14$ \\
\bottomrule
\end{tabular}
\end{center}

The extremal graphs are characterized as follows.
\begin{enumerate}[label=\textup{(\roman*)}]
\item
For $n=10,11,12$, the unique extremal graph for each of
$K_3$ and $K_4$, up to isomorphism, is
$\widehat{\operatorname{Rad}(H)}$, where $H$ is
$W_4$, $C_3\square K_2$, or $W_5$, respectively.

\item
For $n\ge13$ and each $k\in\{3,4\}$, an optimal
$1$-planar graph $G$ with $n$ vertices attains the maximum
of $N(G,K_k)$ if and only if
$G\cong\widehat Q$ for some $Q\in\mathcal G^*$
with $n$ vertices.

\item
For $n\ge15$, an optimal $1$-planar graph $G$ with $n$
vertices attains the maximum of $N(G,K_5)$ if and only if
$G\cong\widehat Q$ for some $Q\in\mathcal G^*$
with $n$ vertices.
\end{enumerate}
\end{theorem}
By Theorem~\ref{thm:cliques}, one sees, for each integer
$n\ge15$, the problems of maximizing $N(G,K_3)$,
$N(G,K_4)$ and $N(G,K_5)$ over optimal $1$-planar graphs
with $n$ vertices have the same extremal graphs.

Finally, for every integer $r\ge1$, we construct a simple
$2$-planar graph $G$ with $5r+3$ vertices and $44r+1$ triangles.
Note that $44r+1>\frac{17(5r+3)-49}{2}.$ Hence, the graph $G$ is a counterexample for Conjecture~\ref{conj2}.

\medskip
\noindent\textbf{Methods and organization.}
Our proofs combine radial graph representations, plane
duality and SPQR trees.
We relate $3$-separators to separating triangles in two
canonical triangulations and obtain the sharp bound by
analyzing the skeletons
of the SPQR tree.
For $K_4$ and $K_5$, we count copies within the bags of
a tree decomposition obtained from the SPQR tree.
To construct the $2$-planar counterexamples, we repeatedly
identify uncrossed facial triangles in copies of a fixed graph.

The organization of this paper is as follows: In Section~\ref{sec:preliminaries}, we provide the definitions and preliminary results.
In Section~\ref{sec:separators}, we prove Theorem~\ref{thm} on the number of $3$-separators
in plane quadrangulations.
In Section~\ref{sec:larger-cliques}, we establish upper bounds for $K_4$ and $K_5$
and complete the proof of Theorem~\ref{thm:cliques}.
Section~\ref{sec:2planar} constructs the counterexample disproving Conjecture~\ref{conj2}.
\section{Preliminaries}
\label{sec:preliminaries}
\subsection{Some notation and preliminaries}
We first give some definitions and preliminary results.
Given a plane graph $G$, its \emph{dual graph} is denoted by $G^*$. Follow \cite[Exercise~10.2.10]{ref03} every plane
quadrangulation is bipartite. An \emph{edge cut} of a graph $G$ is the set of all edges
between $X$ and $V(G)\setminus X$, where $X$ is a nonempty
proper subset of $V(G)$. A \emph{bond} of a graph is a minimal nonempty edge cut.
\begin{lemma}[{\cite{ref03}}]
Every face boundary of a $2$-connected plane graph is a cycle.
\end{lemma}
\begin{lemma}[{\cite{ref03}}]
\label{lem:duality}
Let $G$ be a connected plane graph and let $G^*$ be its dual.
\begin{enumerate}[label=\textup{(\roman*)}]
\item
An edge $e$ of $G$ is a bridge if and only if $e^*$ is
a loop of $G^*$.
\item
A set $C\subseteq E(G)$ is the edge set of a cycle of $G$
if and only if $\{e^*:e\in C\}$ is a bond of $G^*$.
\end{enumerate}
\end{lemma}
\begin{theorem}[{Whitney~\cite{ref03}}]
\label{thm:whitney}
Every simple $3$-connected planar graph has a unique
planar embedding.
\end{theorem}

A triangle $C$ of a graph $G$ is a \emph{separating triangle}
if $V(C)$ is a $3$-separator of $G$.
A cycle of a plane graph is \emph{facial} if it bounds
a face, and \emph{non-facial} otherwise.
The following lemma follows from \cite[Theorem~2.5.1, 
Proposition~5.3.1]{ref48}.

\begin{lemma}[\cite{ref48}]\label{lem:separating-triangle}
Let $T$ be a simple plane triangulation with at least four vertices.
Then $T$ is $3$-connected.
Moreover, a set $X\subseteq V(T)$ with $|X|=3$ is a
$3$-separator of $T$ if and only if $T[X]$ is a non-facial triangle.
\end{lemma}
For a plane triangulation $T$, let $\operatorname{Sep}(T)$ denote the set of separating triangles of $T$, and {set} $\operatorname{sep}(T)=|\operatorname{Sep}(T)|.$
\begin{lemma}[\cite{ref46}]
	\label{lem:radial-connectivity}
	Let $H$ be a connected plane graph. Then $\operatorname{Rad}(H)$ is a $3$-connected quadrangulation
	if and only if $H$ is simple, $2$-connected, and $3$-edge-connected.
\end{lemma}
Alon and Caro~\cite{ref17} and Wood~\cite{ref11}
independently proved the following lemma.
\begin{lemma}[\cite{ref17,ref11}]
\label{lem:planar-K4}
Let $G$ be a { simple} planar graph with at least four vertices.
Then $N(G,K_4)\le |V(G)|-3$.
\end{lemma}
Let $Q$ be a $3$-connected plane quadrangulation with bipartition
$V(Q)=X\mathbin{\cup}Y$.
Then $Q$ is simple (see~\cite{ref46}).
In each face with boundary $xyx'y'x$, where $x,x'\in X$ and
$y,y'\in Y$, draw the diagonals $xx'$ and $yy'$ so that they
cross exactly once.
Let $Q_X$ be the plane graph with vertex set $X$ and edge set
consisting of all diagonals $xx'$.
Define $Q_Y$ similarly on $Y$ with edge set consisting of all
diagonals $yy'$.
Diagonals drawn in different faces are regarded as distinct edges.
\begin{lemma}
\label{lem:radial}
Let $Q_X$ and $Q_Y$ be defined as above. Then 
\begin{enumerate}
\item[\textup{(1)}] $Q_Y=Q_X^*$ and $Q\cong\operatorname{Rad}(Q_X)\cong\operatorname{Rad}(Q_Y)$.
\item[\textup{(2)}] Both $Q_X$ and $Q_Y$ are simple, $2$-connected and $3$-edge-connected.
\item[\textup{(3)}] $|E(Q_X)|=|E(Q_Y)|=|V(Q)|-2$.
\end{enumerate}
\end{lemma}

\begin{proof}
\textup{(1)}
We first describe the faces of $Q_X$.
For each $y\in Y$, let $x_1,\ldots,x_d$ be its neighbours
in cyclic order with $x_{d+1}=x_1$.
Since $Q$ is simple and $3$-connected, these neighbours
are distinct and $d\ge3$.
For each $i$, the face between the consecutive edges
$yx_i$ and $yx_{i+1}$ has $X$-diagonal $x_ix_{i+1}$.
Hence $C_y=x_1x_2\cdots x_dx_1$ is a cycle of $Q_X$. Every face of $Q\cup Q_X$ is a triangle containing
exactly one vertex of $Y$.
Deleting $y$ and its incident edges merges the triangular
faces around $y$ into a face $f_y$ bounded by $C_y$.
Deleting all vertices of $Y$ therefore gives $Q_X$,
whose faces are precisely the distinct faces $f_y$ 
for $y\in Y$.

We now verify duality.
Each face $f_y$ of $Q_X$ contains exactly one vertex of
$Q_Y$, namely $y$.
For a face of $Q$ with boundary $xyx'y'x$, the edge $xx'$
of $Q_X$ is incident with $f_y$ and $f_{y'}$.
The diagonal $yy'$ crosses $xx'$ exactly once and meets
no other edge of $Q_X$, so it is the dual edge of $xx'$.
Hence $Q_Y=Q_X^*$. Moreover, a vertex $x\in X$ lies on the boundary of $f_y$
if and only if $xy\in E(Q)$.
Hence $Q\cong\operatorname{Rad}(Q_X)$ by the definition
of the radial graph.
Interchanging $X$ and $Y$ also gives
$Q\cong\operatorname{Rad}(Q_Y)$.

\smallskip
\noindent
\textup{(2)}
We first prove that $Q_X$ and $Q_Y$ are simple.
Neither graph has loops, since each diagonal joins
distinct vertices of a face of $Q$.
Suppose that $Q_X$ has two parallel edges $e_1,e_2$
with ends $x,x'$.
The edges $e_1$ and $e_2$ form a simple closed curve $\gamma$.
Let $e_1$ be the diagonal of a face of $Q$
with boundary $xyx'y'x$.
The diagonal $yy'$ crosses $e_1$ exactly once.
Since $e_2$ lies in a different face of $Q$, it does not
meet $yy'$.
Thus $y$ and $y'$ lie on opposite sides of $\gamma$.
Since $\gamma$ meets $Q$ only at $x$ and $x'$, there is
no $y$-$y'$ path in $Q-\{x,x'\}$, a contradiction to the
$3$-connectivity of $Q$.
Hence $Q_X$ is simple.
Interchanging $X$ and $Y$ shows that $Q_Y$ is also simple.

We next prove that $Q_X$ is $2$-connected.
Choose three distinct vertices $x, x', z\in X$.
Since $Q-z$ is connected, choose a path
$P=x_0y_1x_1\cdots y_kx_k$ in $Q-z$,
where $x_0=x$, $x_k=x'$, $x_i$ belongs to $X$,
and $y_i$ belong to $Y$.
For each $i=1,\ldots,k$, the vertices $x_{i-1}$ and $x_i$
lie on $C_{y_i}$ and are different from $z$.
Since $C_{y_i}$ is a cycle, choose a path $P_i$
in $C_{y_i}-z$ with ends $x_{i-1}$ and $x_i$.
Replacing the subpath $x_{i-1}y_ix_i$ of $P$ by $P_i$
for each $i$ gives a walk in $Q_X-z$ with ends $x$ and $x'$.
Hence $Q_X-z$ is connected.
Since $|X|\ge3$, $Q_X$ is $2$-connected.
Interchanging $X$ and $Y$ shows that $Q_Y$ is also $2$-connected.

Finally, since $Q_Y=Q_X^*$ is simple,
by Lemma~\ref{lem:duality}, $Q_X$
has no bond of size one or two.
Hence $Q_X$ is $3$-edge-connected.
Similarly, $Q_Y$ is $3$-edge-connected.

\smallskip
\noindent
\textup{(3)}
Since $Q$ is a quadrangulation, $|E(Q)|=2|F(Q)|$.
By Euler's formula, $|F(Q)|=|V(Q)|-2$.
By construction, each of $Q_X$ and $Q_Y$ has exactly
one edge for each face of $Q$.
Hence $|E(Q_X)|=|E(Q_Y)|=|F(Q)|=|V(Q)|-2$.
\end{proof}
\begin{proposition}\label{prop:vertex-face}
A graph $G$ is optimal $1$-planar if and only if
$G\cong\widehat{\operatorname{Rad}(H)}$ for some simple
$2$-connected plane graph $H$ whose dual is simple.
Furthermore, $|V(G)|=|E(H)|+2$, $|E(G)|=4|E(H)|.$
\end{proposition}

\begin{proof}
Let $G$ be an optimal $1$-planar graph. Thus $G\cong\widehat Q$ for some $3$-connected plane quadrangulation $Q$,
where $\widehat Q$ is obtained from $Q$ by adding both diagonals of each face.
By Lemma~\ref{lem:radial}, the graph $H=Q_X$ is simple and
$2$-connected, its dual $H^*=Q_Y$ is simple, and
$Q\cong\operatorname{Rad}(H)$. By Theorem~\ref{thm:whitney}, $Q$ has a unique planar embedding.
Thus
$G\cong\widehat Q\cong\widehat{\operatorname{Rad}(H)}$.

Conversely, let $H$ be a simple $2$-connected plane graph
with simple dual $H^*$. Since $H^*$ is simple, Lemma~\ref{lem:duality}
implies that $H$ is $3$-edge-connected.
Thus $Q=\operatorname{Rad}(H)$ is a $3$-connected plane
quadrangulation by Lemma~\ref{lem:radial-connectivity}.

We first show that $G$ is simple and $1$-planar.
Every face of $Q$ has boundary $ufvgu$, where $uv\in E(H)$
and $f,g$ are the two faces of $H$ incident with $uv$.
In each face $ufvgu$ of $Q$, draw the diagonals $uv$ and $fg$
so that they cross exactly once.
The subgraphs of $G$ induced by $V(H)$ and $F(H)$
are $H$ and $H^*$, respectively, and the edges between
$V(H)$ and $F(H)$ are exactly the edges of $Q$.
Since $Q$, $H$ and $H^*$ are simple, $G$ is simple.
In this drawing, each edge of $H$ crosses its dual edge
in $H^*$ exactly once.
Hence $G$ is $1$-planar.

By the definition of the radial graph and Euler's formula, $|V(G)|=|V(Q)|=|V(H)|+|F(H)|=|E(H)|+2.$
Since $Q$ is a quadrangulation,
$|E(Q)|=2|V(Q)|-4=2|E(H)|$.
Adding the diagonals introduces the edges of $H$ and $H^*$.
Since $|E(H^*)|=|E(H)|$, we obtain
\[
|E(G)|=|E(Q)|+|E(H)|+|E(H^*)|
=4|E(H)|=4|V(G)|-8.
\]
Therefore $G$ is optimal $1$-planar.
\end{proof}

\subsection{SPQR tree and some basic properties}
Recall the construction of the SPQR tree for a connected graph \(G\) with at least three edges and no cut vertex;
see \cite{Angelini2014,ref06,ref08}.
In this construction, $G$ may be a graph in which two vertices are joined by at least three parallel edges. A pair of distinct vertices \(u, v \in V(G)\) is called a \textit{separation pair} if \(G - \{u, v\}\) is disconnected. A \textit{split pair} is defined as either a separation pair or a pair of adjacent vertices.

Let $\{u,v\}$ be a split pair of $G$. Suppose that
$G=G_1\cup G_2$, $E(G_1)\cap E(G_2)=\emptyset$,
$V(G_1)\cap V(G_2)=\{u,v\}$, and $|E(G_i)|\ge 2$
for $i=1,2$. A \emph{split} at $\{u,v\}$ replaces
$G$ by $G'_1=G_1+e_1$ and $G'_2=G_2+e_2$, where
$e_1$ and $e_2$ are new virtual edges with ends $u$
and $v$. We call $e_1$ and $e_2$ \emph{twin edges}.
The ends of a virtual edge are called its \emph{poles}.
Apply the same operation to the resulting graphs
until no further split is possible.
By \cite{ref08}, the resulting graphs are triangles,
graphs with two vertices joined by three parallel
edges, or simple 3-connected graphs. We call them
\emph{split components}.

Two graphs are called \emph{adjacent} if they contain
a pair of twin virtual edges. If two adjacent graphs
are both cycles or both have two vertices, merge them
by deleting the twin virtual edges and identifying
the corresponding ends. Repeat until no such pair
remains. Throughout the construction, vertices retain
their original labels, and surviving virtual edges
retain their twin relations.
The obtained graphs are called \emph{skeletons}.
An $S$-skeleton is a cycle of length at least three.
A $P$-skeleton has two vertices and at least three
parallel edges. An $R$-skeleton is a simple
3-connected graph. Each edge of graph
$G$ belongs to exactly one skeleton and is called
a \emph{real edge}.

We now construct a graph $\mathcal T$. Take one node
for each skeleton, and call it an S-, P-, or R-node
according to the type of the skeleton.
For each pair of twin virtual edges, join the nodes
whose skeletons contain the two edges.
For each real edge $e\in E(G)$, add a new node $q_e$
and join it to the node whose skeleton contains $e$.
We call $q_e$ a \emph{Q-node}.

By \cite{ref06,ref08}, $\mathcal T$ is a tree whose
leaves are exactly the Q-nodes.
For each S-, P-, or R-node, every real edge of its
skeleton corresponds to a tree edge joining it to
a Q-node, and every virtual edge corresponds to a
tree edge joining it to another S-, P-, or R-node.
No two S-nodes are adjacent, and no two P-nodes are
adjacent. We call $\mathcal T$ the \emph{SPQR tree}
of $G$. To recover $G$, take disjoint copies of
the skeletons, retaining their vertex labels. For each pair of twin
virtual edges, identify the ends with the same label and delete
the two virtual edges.

Under the SPQR-tree convention of Angelini et al.~\cite{Angelini2014}, the skeleton of the Q-node \(q_e\) contains the real edge $e$ together with a parallel virtual edge. This virtual edge has a twin virtual edge inside the adjacent skeleton. Our convention is obtained by replacing this twin with $e$ while keeping \(q_e\) as a leaf node.

\begin{lemma}[{\cite{Angelini2014}}]
\label{lem:spqr-uniqueness}
Let $G$ be a connected planar graph without cut vertex and with at least
three edges. Then the SPQR tree of $G$ is unique up to isomorphism,
regardless of the planar embedding chosen for $G$.
\end{lemma}

\begin{lemma}[{\cite{Angelini2014}}]
\label{lem:spqr-duality}
Let $G$ be a $2$-connected planar graph with a fixed planar embedding, and let $G^*$ be its dual. Let $\mathcal{T}$ and $\mathcal{T}^*$ denote the SPQR trees of $G$ and $G^*$, respectively. Then there exists a tree isomorphism $\phi\colon\mathcal{T}\to\mathcal{T}^*$ such that $\phi(q_e)=q_{e^*}$ for all $e\in E(G)$, where $e^*$ is the dual edge of $e$. This isomorphism swaps S-nodes and P-nodes, and sends each R-node $x\in\mathcal{T}$ to an R-node $\phi(x)\in\mathcal{T}^*$ whose skeleton is dual to that of $x$.
\end{lemma}
Let $H$ be a simple $2$-connected plane graph with a simple dual $H^*$, and let $\mathcal{T}$ denote its SPQR tree. Let $\mathcal{T}^\circ$ be the tree obtained from $\mathcal{T}$ by removing all $Q$-nodes; we refer to $\mathcal{T}^\circ$ as the \emph{core} of $\mathcal{T}$. Write $\mathcal{R}$ for the collection of $R$-nodes of $\mathcal{T}$. For each $x\in V(\mathcal{T}^\circ)$, let $J_x$ stand for the skeleton of $x$, equipped with the embedding inherited from $H$, and let $q_x$ be the number of $Q$-node neighbors of $x$ in $\mathcal{T}$. Define
\(
d_x:=\deg_{\mathcal{T}}(x)=|E(J_x)|,\ \ c_x:=\deg_{\mathcal{T}^\circ}(x)=d_x-q_x.
\)
Now set
\[
\begin{aligned}
U&:=\{x\in V(\mathcal{T}^\circ):
      x\text{ is an $S$- or $P$-node and }q_x=1\},\\
W&:=\{x\in V(\mathcal{T}^\circ):
      x\text{ is an $S$- or $P$-node and }q_x=0\}.
\end{aligned}
\]

\begin{lemma}\label{lem:core-leaves}
Let $H$ be a simple $2$-connected plane graph with simple dual $H^*$, and let $\mathcal{T}^\circ$ denote the core of its SPQR tree.
Write $Q=\operatorname{Rad}(H)$ and $n=|V(Q)|$. Then
every vertex of $\mathcal{T}^\circ$ of degree at most one is an $R$-node.
Moreover,
\begin{equation}\label{eq:spqr-size}
n=4+\sum_{x\in V(\mathcal{T}^\circ)}(d_x-2)
  \ge 4+4|\mathcal R|+|U|+|W|.
\end{equation}
\end{lemma}

\begin{proof}
The edges of a $P$-skeleton are parallel. Since both $H$ and $H^*$ are simple, every $P$-skeleton in either SPQR tree contains at most one real edge. Consequently, each $P$-node in either tree has at most one $Q$-node neighbor. By Lemma~\ref{lem:spqr-duality}, every $S$-node of $\mathcal{T}$ maps to a $P$-node in the SPQR tree of $H^*$ with the same number of $Q$-node neighbors. It follows that $q_x\le 1$ for all $S$-nodes and $P$-nodes $x$ of $\mathcal{T}$.

Every $S$-skeleton and every $P$-skeleton has at least three edges. Hence $c_x=d_x-1\ge 2$ for $x\in U$ and $c_x=d_x\ge 3$ for $x\in W$. Since $U\cup W$ exhausts all $S$-nodes and $P$-nodes of $\mathcal{T}$, any vertex of $\mathcal{T}^\circ$ of degree at most one must be an $R$-node.

Let $m=|E(H)|$. By Euler's formula, $n=|V(H)|+|F(H)|=m+2$. As $\mathcal{T}$ is a tree,
\(
\sum_{x\in V(\mathcal{T})}\bigl(\deg_{\mathcal{T}}(x)-2\bigr)=-2.
\)
Each of the $m$ $Q$-nodes contributes $-1$ to this sum, so
\[
\sum_{x\in V(\mathcal{T}^\circ)}(d_x-2)=m-2=n-4.
\]
Every $R$-skeleton has at least six edges, while every $S$-skeleton or $P$-skeleton has at least three edges. Therefore $d_x-2\ge 4$ for $x\in\mathcal{R}$ and $d_x-2\ge 1$ for $x\in U\cup W$. This proves the inequality in \eqref{eq:spqr-size}.
\end{proof}

\section{Proof of Theorem~\ref{thm}}\label{sec:separators}
\subsection{\texorpdfstring{$3$}{3}-separators in plane quadrangulations}

Let $Q$ be a $3$-connected plane quadrangulation with bipartition $V(Q)=X\cup Y$.
Its two \emph{canonical triangulations} are $T_X:=Q\cup Q_X$ and $T_Y:=Q\cup Q_Y$.
Explicitly, $T_X$ is formed by adding all $X$-diagonals of $Q$, and $T_Y$ by adding all $Y$-diagonals.
By Lemma~\ref{lem:radial}, both are simple plane triangulations.

\begin{lemma}\label{lem:separators}
Let $S\subseteq V(Q)$. Then $S$ is a $3$-separator of $Q$ if and only if $S$ induces a separating triangle in $T_X$ or in $T_Y$. Moreover,
\(
s_3(Q)=\operatorname{sep}(T_X)+\operatorname{sep}(T_Y).
\)
\end{lemma}

\begin{proof}
Let $S$ be a $3$-separator of $Q$ with $|S\cap X|\ge 2$, so $|S\cap Y|\le 1$.
Consider any $X$-diagonal $xx'$ in $T_X-S$. The face of $Q$ containing this diagonal has boundary $xyx'y'x$ for some $y,y'\in Y$.
At least one of $y,y'$ lies outside $S$, so $x$ and $x'$ are connected by a path of length two in $Q-S$.
It follows that none of the added diagonals connect distinct components of $Q-S$.
Hence $T_X-S$ is disconnected, and $S$ is a $3$-separator of $T_X$.
By Lemma~\ref{lem:separating-triangle}, $S$ induces a separating triangle in $T_X$.
Symmetrically, if $|S\cap Y|\ge 2$, then $S$ induces a separating triangle in $T_Y$.

Conversely, suppose $S$ induces a separating triangle in $T_X$.
Since $T_X-S$ is disconnected and $Q-S$ is a spanning subgraph of $T_X-S$, $Q-S$ is also disconnected.
Thus $S$ is a $3$-separator of $Q$.
Because every edge of $T_X$ added to $Q$ is an $X$-diagonal, $T_X$ contains no edge with both ends in $Y$, which forces $|S\cap X|\ge 2$.
The same reasoning applies to $T_Y$.

Each $3$-separator of $Q$ belongs to exactly one of the two triangulations, determined by which partition class contains at least two of its vertices.
We therefore conclude
\(
s_3(Q)=\operatorname{sep}(T_X)+\operatorname{sep}(T_Y).
\)
\end{proof}

\begin{corollary}\label{cor:triangles}
Let $Q$ be a $3$-connected plane quadrangulation on $n$ vertices. Then
\(
N(\widehat Q,K_3)=4n-8+s_3(Q).
\)
\end{corollary}

\begin{proof}
Every triangle of $\widehat Q$ lies in precisely one of $T_X$ and $T_Y$.
Therefore, $N(\widehat Q,K_3)=N(T_X,K_3)+N(T_Y,K_3)$.
Both $T_X$ and $T_Y$ are simple plane triangulations on $n$ vertices. By Euler's formula, each contains $2n-4$ facial triangles.
By Lemma~\ref{lem:separating-triangle}, the non-facial triangles of each triangulation are exactly its separating triangles.
Hence $N(T_X,K_3)=2n-4+\operatorname{sep}(T_X)$ and $N(T_Y,K_3)=2n-4+\operatorname{sep}(T_Y)$.
Combining this with Lemma~\ref{lem:separators}, we obtain
\(
N(\widehat Q,K_3)=4n-8+s_3(Q).
\)
\end{proof}

Let $G$ be a simple $2$-connected plane graph, and let $uv\in E(G)$.
The edge $uv$ itself is called the \emph{trivial $uv$-bridge}.
For each component $K$ of $G-\{u,v\}$, the subgraph with vertex set $V(K)\cup\{u,v\}$ and edge set $E(K)$ together with all edges joining $K$ to $u$ or $v$ is called a \emph{non-trivial $uv$-bridge}.

\begin{lemma}\label{lem:common-faces}
Let $G$ be a simple $2$-connected plane graph and let $uv\in E(G)$.
If $G-\{u,v\}$ has $k$ components, then precisely $k+1$ face boundaries contain both $u$ and $v$, and exactly two of these face boundaries contain the edge $uv$.
\end{lemma}

\begin{proof}
Let $B_1,\dots,B_k$ denote the non-trivial $uv$-bridges.
Call a face \emph{common} if its boundary contains both $u$ and $v$.
Fix $i\in\{1,\dots,k\}$.
The edge $uv$ lies in some common face of $B_i$; we prove that this common face is unique.

Suppose, for contradiction, that $B_i$ admits two distinct common faces $f_1$ and $f_2$.
Since $u$ belongs to two distinct face boundaries, $u$ has at least two neighbors in $B_i$.
Choose consecutive edges $ua,ub$ around $u$ such that $f_1$ lies between them.
As $B_i-\{u,v\}$ is connected, there exists a path $P\subseteq B_i-\{u,v\}$ with ends $a$ and $b$.
Then $C=ua\cup P\cup ub$ is a cycle avoiding $v$.
The faces $f_1$ and $f_2$ lie on opposite sides of $C$ near $u$.
Because each face is connected and disjoint from $C$, $f_1$ and $f_2$ lie on different sides of $C$.
But both face boundaries contain $v$, which forces $v\in V(C)$, a contradiction.
Hence each $B_i$ has exactly one common face.

Let $G_0$ consist of the single edge $uv$.
For $i=1,\dots,k$, set $G_i:=G_{i-1}\cup B_i$, so $G_k=G$.
For each $1\le j\le k$, $B_j-\{u,v\}$ is connected and has neighbors at both $u$ and $v$; this follows from the $2$-connectivity of $G$.
Fix $i\in\{1,\dots,k\}$, and let $f$ be the unique common face of $B_i$.
For all $j<i$, $B_j-\{u,v\}$ lies inside a face of $B_i$ whose boundary contains both $u$ and $v$.
By uniqueness, this face must be $f$, and $uv$ also belongs to $f$.
Thus $G_{i-1}$ is contained in $f$, with $u$ and $v$ on its boundary.
Symmetrically, $B_i-\{u,v\}$ lies inside some common face $g$ of $G_{i-1}$.

In $G_i$, all faces of $B_i$ other than $f$ and all faces of $G_{i-1}$ other than $g$ are preserved.
Since $G_{i-1}$ and $B_i$ share no edges, $G_i$ inherits $|F(G_{i-1})|+|F(B_i)|-2$ such unchanged faces.
Both $G_{i-1}$ and $B_i$ are connected, $V(G_{i-1})\cap V(B_i)=\{u,v\}$, and $E(G_{i-1})\cap E(B_i)=\emptyset$.
By Euler's formula, $|F(G_i)|=|F(G_{i-1})|+|F(B_i)|$.
It follows that $G_i$ contains exactly two new faces.

For every $1\le j\le i$, $B_j-u$ is connected and contains $v$.
Hence $G_i-u$ is connected, and similarly $G_i-v$ is connected.
Consequently, every face boundary of $G_i$ passes through each of $u$ and $v$ at most once.
Each new face contains edges from both $G_{i-1}$ and $B_i$.
Since these subgraphs intersect only at $u$ and $v$, its boundary contains at least one of $u$ or $v$.
If a new face boundary contained only $u$, its boundary walk would visit $u$ at least twice, which is impossible.
Therefore the two new faces are common faces.
The common face $g$ of $G_{i-1}$ splits into exactly two common faces in $G_i$.
Thus $G_i$ has one more common face than $G_{i-1}$.

Since $G_0$ has one common face, $G_k=G$ has precisely $k+1$ common faces.
As $G$ is $2$-connected, the edge $uv$ is incident to two distinct faces, so exactly two of these common faces contain $uv$.
\end{proof}

For $uv\in E(G)$, let $c_{uv}$ denote the number of components of $G-\{u,v\}$.
For a simple $2$-connected plane graph $G$, define
\[
\eta(G):=
\bigl|\{(uv,f):uv\in E(G),\ f\in F(G),\
u,v\in V(\partial f),\ uv\notin E(\partial f)\}\bigr|.
\]
By Lemma~\ref{lem:common-faces}, $\eta(G)=\sum_{uv\in E(G)}(c_{uv}-1)$.

\begin{lemma}\label{lem:four-term}
Let $Q$ be a $3$-connected plane quadrangulation with bipartition $V(Q)=X\mathbin{\cup}Y$. Then
\(
s_3(Q)=N(Q_X,K_3)+N(Q_Y,K_3)+\eta(Q_X)+\eta(Q_Y).
\)
The four terms count the $3$-separators of types $XXX$, $YYY$, $XXY$, and $XYY$, respectively.
\end{lemma}

\begin{proof}
By Lemma~\ref{lem:separators}, $s_3(Q)=\operatorname{sep}(T_X)+\operatorname{sep}(T_Y)$.
We first enumerate the separating triangles of $T_X$. Since $Y$ is an independent set in $T_X$, every triangle of $T_X$ is of type $XXX$ or $XXY$.
The $XXX$-triangles in $T_X$ are precisely the triangles of $Q_X$.
Every face of $T_X$ has type $XXY$, so every $XXX$-triangle is non-facial and therefore separating by Lemma~\ref{lem:separating-triangle}.
Thus $T_X$ contains $N(Q_X,K_3)$ separating triangles of type $XXX$.

Now take $uv\in E(Q_X)$, and let $f$ be a face of $Q_X$ whose boundary contains both $u$ and $v$. Let $y_f$ denote the vertex of $Y$ lying in $f$.
Then $uvy_f$ forms an $XXY$-triangle in $T_X$.
Conversely, the two vertices in $X$ of any $XXY$-triangle determine the edge $uv$, and its unique vertex in $Y$ determines the corresponding face $f$.

Let $f_1,f_2$ be the two distinct faces of $Q_X$ incident with $uv$.
The face of $Q$ associated with $uv$ has boundary $uy_{f_1}vy_{f_2}u$.
Adding the diagonal $uv$ yields the two facial triangles of $T_X$ containing $uv$, namely $uvy_{f_1}$ and $uvy_{f_2}$.
Hence $uvy_f$ is facial in $T_X$ if and only if $uv\in E(\partial f)$.
By Lemma~\ref{lem:separating-triangle}, it is separating if and only if $uv\notin E(\partial f)$.
It follows that $T_X$ has $\eta(Q_X)$ separating triangles of type $XXY$, so
\[
\operatorname{sep}(T_X)=N(Q_X,K_3)+\eta(Q_X).
\]
Symmetrically, $\operatorname{sep}(T_Y)=N(Q_Y,K_3)+\eta(Q_Y)$.
Combining these identities with Lemma~\ref{lem:separators} yields the claimed equality.
\end{proof}

For an internal node $x$ of the SPQR tree of a graph $G$ and an edge $e$ of the skeleton of $x$, we define a subgraph $G_{x,e}$ of $G$ as follows.
If $e$ is real, set $G_{x,e}:=e$.
If $e$ is virtual, remove the tree edge corresponding to its twin pair and take the component not containing $x$.
Recompose the skeletons within this component and delete the virtual edge paired with $e$; the resulting subgraph is denoted $G_{x,e}$.
For fixed $x$, each edge of $G$ belongs to exactly one subgraph $G_{x,e}$. Moreover, $G_{x,e}$ intersects the other such subgraphs only at the ends of $e$.

\begin{lemma}\label{lem:mixed-spqr}
Let $H$ be a simple $2$-connected plane graph for which the dual $H^*$ is also simple, and let $\mathcal{T}$ be its SPQR tree.
Set $m:=|E(H)|$, $r:=\eta(H)+\eta(H^*)$, and
\(
C:=\sum_{\rho\in\mathcal R}(d_\rho-2)+\sum_{x\in W}(d_x-2).
\)
Then $r=m-2-C$.
\end{lemma}

\begin{proof}
For each internal node $x$ of $\mathcal{T}$, let $J_x$ denote its skeleton.
First suppose $x$ is an $S$-node or an $R$-node, and let $uv$ be a real edge of $J_x$.
Then $J_x-\{u,v\}$ is nonempty and connected.
For each virtual edge $e=ab$ of $J_x$, the subgraph $H_{x,e}$ meets $J_x$ precisely at $a$ and $b$.
By the split construction, $H_{x,e}+ab$ is $2$-connected, so $H_{x,e}$, $H_{x,e}-a$, and $H_{x,e}-b$ are connected.
Since $J_x$ is simple, at most one of $a,b$ lies in $\{u,v\}$.
It follows that $H_{x,e}-\{u,v\}$ is connected and intersects $J_x-\{u,v\}$.
A path in $J_x-\{u,v\}$ lifts to a walk in $H-\{u,v\}$ by replacing each virtual edge $e$ along the path with a path in $H_{x,e}-\{u,v\}$.
Thus all vertices of $J_x-\{u,v\}$ belong to a single component of $H-\{u,v\}$.
Every remaining vertex of $H-\{u,v\}$ lies in some connected subgraph $H_{x,e}-\{u,v\}$ that intersects this component.
Hence $H-\{u,v\}$ is connected, so $c_{uv}=1$.

Now let $x$ be a $P$-node possessing a real edge $uv$.
As $H$ is simple, the remaining $d_x-1$ edges of $J_x$ are virtual; label them $e_1,\dots,e_{d_x-1}$.
For each $i$, let $y_i$ be the neighbor of $x$ whose skeleton contains the twin of $e_i$.
No two $P$-nodes are adjacent, so each $y_i$ is an $S$-node or an $R$-node.
Since $J_{y_i}-\{u,v\}$ is nonempty and connected, the argument above implies that $H_{x,e_i}-\{u,v\}$ is nonempty and connected.
These subgraphs are pairwise disjoint, and their union equals $H-\{u,v\}$.
Consequently, $c_{uv}=d_x-1$.

Recall that $\eta(H)=\sum_{uv\in E(H)}(c_{uv}-1)$.
An edge contributes positively to $\eta(H)$ if and only if it is a real edge of a $P$-skeleton.
Every $P$-node in $U$ has exactly one real edge, so
\[
\eta(H)=\sum_{\substack{x\in U\\ x\text{ is a }P\text{-node}}}(d_x-2).
\]
Applying the same reasoning to $H^*$ and invoking Lemma~\ref{lem:spqr-duality}, we obtain
\[
\eta(H^*)=\sum_{\substack{x\in U\\ x\text{ is an }S\text{-node}}}(d_x-2).
\]
Therefore $r=\sum_{x\in U}(d_x-2)$.

The $Q$-nodes are exactly the leaves of $\mathcal{T}$, so $\mathcal{T}$ has $m$ leaves.
By the handshaking lemma for trees,
\(
\sum_{x\in V(\mathcal{T}^\circ)} (d_x-2)=m-2.
\)
The sets $U$, $W$, and $\mathcal{R}$ partition the internal nodes of the core. Hence
\(
m-2
=\sum_{x\in U}(d_x-2)
 +\sum_{x\in W}(d_x-2)
 +\sum_{\rho\in\mathcal R}(d_\rho-2)
=r+C.
\)
Rearranging yields $r=m-2-C$.
\end{proof}

\subsection{Triangles in SPQR skeletons}

\begin{lemma}\label{lem:cycle-projection}
Let $G$ be a $2$-connected graph with SPQR tree $\mathcal{T}$, and let $x$ be an internal node of $\mathcal{T}$.
Let $C$ be a cycle of $G$ that intersects at least two distinct subgraphs $G_{x,e}$.
Define $P_x(C)$ as the subgraph of the skeleton of $x$ whose edge set consists of all edges $e$ satisfying $E(C)\cap E(G_{x,e})\neq\emptyset$, and whose vertex set is the set of ends of these edges.
Then $P_x(C)$ is a cycle in the skeleton of $x$.
\end{lemma}
\begin{proof}
Let $e=uv$ be an edge of $P_x(C)$.
Then $C$ contains at least one edge from $G_{x,e}$.
Because $C$ meets more than one subgraph $G_{x,e}$, it also contains edges outside $G_{x,e}$.
The edges of $C$ in $G_{x,e}$ meet the other edges of $C$
only at $u$ and $v$. Hence the edges of $C$ in $G_{x,e}$ form a $u$--$v$ path.
For each edge $e\in E(P_x(C))$, replace this $u$--$v$ subpath of $C$ within $G_{x,e}$ by the skeleton edge $e$.
The resulting graph is exactly $P_x(C)$.
Since every path is substituted by an edge with identical ends, $P_x(C)$ is a cycle in the skeleton of $x$.
If $x$ is an $S$-node, then $P_x(C)$ coincides with the entire skeleton of $x$; if $x$ is a $P$-node, then $P_x(C)$ consists of a pair of parallel edges.
\end{proof}

Let $H$ be a simple $2$-connected plane graph whose dual $H^*$ is also simple, and let $\mathcal{T}$ denote its SPQR tree.
Let $x$ be an $S$-node or $R$-node of $\mathcal{T}$, and let $e=uv$ be an edge of its skeleton.
We call $e$ \emph{$H$-active} if either
\begin{enumerate}
    \item[(i)] the tree edge associated with $e$ connects $x$ to a Q-node; or
    \item[(ii)] the tree edge associated with $e$ connects $x$ to a $P$-node $y\in U$.
\end{enumerate}
Every $H$-active edge $e=uv$ corresponds to a real edge of $H$.
In case (i), $e$ itself is a real edge of $H$.
In case (ii), since $y\in U$, the skeleton of $y$ contains exactly one real edge $f$.
All edges of this $P$-skeleton share the ends $u$ and $v$, so $e$ corresponds to the real edge $f=uv$ of $H$.

A triangle in the skeleton of an $S$-node or $R$-node of $\mathcal{T}$ is called \emph{fully $H$-active} if all three of its edges are $H$-active.
We define $H^*$-active edges and fully $H^*$-active triangles analogously with respect to the SPQR tree of $H^*$.

\begin{lemma}\label{lem:triangle-home}
For every triangle $C$ of $H$, there exists a unique pair $(x,D)$ such that $x\in V(\mathcal{T})$ and $D$ is a fully $H$-active triangle in the skeleton of $x$ satisfying $V(D)=V(C)$.
We refer to $x$ as the \emph{home} of $C$.
The node $x$ is either an $R$-node or an $S$-node of degree three in $\mathcal{T}$.
Conversely, if $D$ is a fully $H$-active triangle in the skeleton of some node $x$, then $H[V(D)]$ is a triangle whose home is $x$.
The same conclusion holds for $H^*$ using its dual skeletons.
\end{lemma}
\begin{proof}
Let $C$ be a triangle of $H$, and let $q_1,q_2,q_3$ be the Q-nodes corresponding to its three edges.
Consider the minimal subtree of $\mathcal{T}$ containing $q_1,q_2,q_3$.
Its leaves are exactly $q_1,q_2,q_3$, so this subtree contains a unique node $x$ of degree three.
Thus $q_1,q_2,q_3$ lie in three distinct connected components of $\mathcal{T}-x$.
The three tree edges incident to $x$ within this subtree correspond to three edges in the skeleton of $x$.
By Lemma~\ref{lem:cycle-projection}, these three skeleton edges form a cycle and hence a triangle, which we denote by $D$.
Since a $P$-skeleton has only two vertices, $x$ cannot be a $P$-node; therefore $x$ is an $R$-node or an $S$-node.
If $x$ is an $S$-node, its skeleton is precisely $D$, which forces $\deg_{\mathcal{T}}(x)=3$.

If $x$ is the only internal node of this minimal subtree, then $q_1,q_2,q_3$ are all adjacent to $x$, and the three edges of $D$ are $H$-active.
Suppose instead that the minimal subtree contains another internal node $y\neq x$.
Then $\deg(y)=2$ within this subtree.
The two tree edges incident to $y$ correspond to two edges of the skeleton of $y$.
By Lemma~\ref{lem:cycle-projection}, these two edges form a $2$-cycle, so $y$ must be a $P$-node.

No two $P$-nodes are adjacent in $\mathcal{T}$.
It follows that the path from $x$ to each $q_i$ is either the direct edge $xq_i$, or a path $xy_iq_i$ where $y_i$ is a $P$-node.
In the latter case, simplicity of $H$ implies that $q_i$ is the unique Q-leaf adjacent to $y_i$, so $y_i\in U$.
By the definition of $H$-active edges, all three edges of $D$ are $H$-active, and their associated real edges are exactly the three edges of $C$.
Because each edge of $D$ has the same ends as its corresponding edge in $C$, we obtain $V(D)=V(C)$.

Now suppose $D'$ is another fully $H$-active triangle in the skeleton of a node $z$ with $V(D')=V(C)$.
As $H$ is simple, the real edges associated with $D'$ coincide with the edges of $C$.
From the definition of $H$-active edges, the Q-nodes $q_1,q_2,q_3$ lie in three distinct components of $\mathcal{T}-z$.
Thus $z$ is the unique degree-three node in the minimal subtree containing $q_1,q_2,q_3$, and so $z=x$.
For each $i\in\{1,2,3\}$, the component of $\mathcal{T}-x$ containing $q_i$ attaches to $x$ via a unique tree edge, whose corresponding skeleton edge belongs to both $D$ and $D'$.
Hence $D$ and $D'$ share the same three edges, giving $D'=D$.

Conversely, let $D$ be a fully $H$-active triangle in the skeleton of $x$, with vertices $u,v,w$.
Each edge of $D$ has the same ends as its associated real edge in $H$, so $uv,vw,wu\in E(H)$.
Therefore $H[V(D)]$ is a triangle, and by uniqueness its home is $x$.

By Lemma~\ref{lem:spqr-duality}, the same reasoning applies to $H^*$.
Consequently, every triangle of $H^*$ admits a unique home, which is either an $R$-node or an $S$-node of degree three in the SPQR tree of $H^*$.
\end{proof}

For a plane graph $H$, define $\delta(H):=|E(H)|-2-N(H,K_3)$ and $\Delta(H):=2|E(H)|-4-N(H,K_3)-N(H^*,K_3)$.
Since $|E(H)|=|E(H^*)|$, it follows that $\Delta(H)=\delta(H)+\delta(H^*)$.

Now suppose $H$ is simple and $2$-connected, and its dual $H^*$ is also simple.
Let $\mathcal{T}$ denote the SPQR tree of $H$.
For each $x\in V(\mathcal{T}^\circ)$, let $\ell_x$ be the sum of the number of triangles in $J_x$ that are not fully $H$-active and the number of triangles in $J_x^*$ that are not fully $H^*$-active.
For each $S$-node or $P$-node $x$, let $z_x$ denote the number of triangles of $H$ and $H^*$ with home $x$.
By Lemma~\ref{lem:triangle-home}, $z_x\in\{0,1\}$.

\begin{lemma}\label{lem:identity}
Let $Q$ be a $3$-connected plane quadrangulation with bipartition $V(Q)=X\mathbin{\cup}Y$, and let $\mathcal{T}$ be the SPQR tree of $Q_X$.
Set $n:=|V(Q)|$ and $\mu(Q):=2n-8-s_3(Q)$.
Then
\begin{equation}\label{eq:potential}
    \mu(Q)
    =\sum_{\rho\in\mathcal R}\bigl(\Delta(J_\rho)+\ell_\rho\bigr)
    +\sum_{x\in U}(c_x-1-z_x)
    +\sum_{x\in W}(2c_x-4-z_x).
\end{equation}
Moreover, $c_x-1-z_x\ge0$ for all $x\in U$, and
$2c_x-4-z_x\ge1$ for all $x\in W$.
\end{lemma}
\begin{proof}
By Lemma~\ref{lem:triangle-home}, every triangle of $Q_X$ or $Q_Y$ has a unique home.
Therefore
\begin{equation}\label{eq:p-home}
N(Q_X,K_3)+N(Q_Y,K_3)
=\sum_{\rho\in\mathcal R}
\bigl(N(J_\rho,K_3)+N(J_\rho^*,K_3)-\ell_\rho\bigr)
+\sum_{x\in U\cup W}z_x.
\end{equation}

The Q-nodes are precisely the leaves of $\mathcal{T}$.
Let $m:=|E(Q_X)|$; there are exactly $m$ Q-nodes.
Each node in $U$ is adjacent to a Q-node, while nodes in $W$ have no Q-node neighbors.
Hence
\begin{equation}\label{eq:m-q}
m=\sum_{\rho\in\mathcal R}q_\rho+|U|.
\end{equation}
By Lemmas~\ref{lem:four-term} and~\ref{lem:mixed-spqr}, $s_3(Q)=N(Q_X,K_3)+N(Q_Y,K_3)+m-2-C$,
where $C$ is the quantity defined in Lemma~\ref{lem:mixed-spqr}.
Since $n=m+2$, we have $\mu(Q)=m+C-N(Q_X,K_3)-N(Q_Y,K_3)-2$.
Together with \eqref{eq:p-home}, \eqref{eq:m-q}, and the definition of $C$, we have 
\begin{equation}\label{eq5}
\mu(Q)
=\sum_{\rho\in\mathcal R}
\bigl(q_\rho+d_\rho-2-N(J_\rho,K_3)-N(J_\rho^*,K_3)+\ell_\rho\bigr)
+\sum_{x\in U}(1-z_x)
+\sum_{x\in W}(d_x-2-z_x)-2.
\end{equation}
Recall $q_\rho=d_\rho-c_\rho$.
The summand indexed by $\rho$ then simplifies to $\Delta(J_\rho)+\ell_\rho-(c_\rho-2)$.
For $x\in U$, $1-z_x=(c_x-1-z_x)-(c_x-2)$.
For $x\in W$, $d_x=c_x$, so
$d_x-2-z_x=(2c_x-4-z_x)-(c_x-2)$.
Substituting these identities into \eqref{eq5} yields
\[
\mu(Q)
=\sum_{\rho\in\mathcal R}(\Delta(J_\rho)+\ell_\rho)
+\sum_{x\in U}(c_x-1-z_x)
+\sum_{x\in W}(2c_x-4-z_x)
-\sum_{x\in V(\mathcal T^\circ)}(c_x-2)-2.
\]
Since $\mathcal{T}^\circ$ is a tree, the handshaking lemma for trees gives
\[
\sum_{x\in V(\mathcal T^\circ)}(c_x-2)
=2|E(\mathcal T^\circ)|-2|V(\mathcal T^\circ)|=-2.
\]
This establishes the desired identity.

Now verify the lower bounds.
If $x\in U$, then $c_x\ge2$ and $z_x\le1$, so $c_x-1-z_x\ge0$.
If $x\in W$, then $c_x\ge3$ and $z_x\le1$, so $2c_x-4-z_x\ge1$.
\end{proof}

\subsection{Bounds for R-skeletons}
\label{sec:skeleton-bounds}
Let $G$ be a simple $3$-connected plane graph. Then its dual $G^*$ is also simple and $3$-connected.
A triangle of $G$ is \emph{separating} if both of the regions it bounds contain vertices of $G$.
Indeed, after deleting the vertices of a separating triangle, each resulting component has neighbors on all three vertices of the triangle; planarity ensures that at most one such component lies on each side of the triangle.

Let $C$ be a separating triangle of $G$.
Let $G_1$ and $G_2$ be the subgraphs obtained from $G$ by retaining, respectively, all vertices lying on and outside $C$, and all vertices lying on and inside $C$.
Then $G=G_1\cup G_2$ and $G_1\cap G_2=C$.
We say that $G$ is \emph{decomposed along $C$} into $G_1$ and $G_2$.
Iterate this decomposition procedure until no subgraph in the collection contains a separating triangle.
We call the resulting graphs the \emph{atoms} of $G$.
For each atom $A$, write $v_A=|V(A)|$, and let $q_A$ denote the number of non-triangular faces of $A$.

Let $G$ be a plane graph, and let $f$ be a triangular face
with boundary $abc$. A \emph{stacking in $f$} consists of inserting
a new vertex $v$ into $f$ and adding the edges $va,vb,vc$,
with pairwise disjoint interiors contained in $f$. 
\begin{lemma}\label{lem:triangle-atoms}
Every atom $A$ is simple and $3$-connected. Moreover,
\(
\delta(G)=\sum_A (v_A-4+q_A),
\)
where the sum runs over all atoms of $G$.
If every atom except a distinguished atom $A_0$ is isomorphic to $K_4$, then $G$ can be constructed from $A_0$ by a sequence of stackings.
\end{lemma}

\begin{proof}
Suppose a simple $3$-connected plane graph $G$ is decomposed along a separating triangle $C$ into $G_1$ and $G_2$.
Both $G_1$ and $G_2$ are simple and have at least four vertices.
Take $i\in\{1,2\}$ and any subset $S\subseteq V(G_i)$ with $|S|\le 2$.
Since $C-S$ is nonempty and connected, every connected component of $G_i-S$ intersects $C-S$; otherwise $G-S$ would be disconnected.
It follows that $G_i-S$ is connected.
Hence $G_1$ and $G_2$ are $3$-connected, and therefore every atom is simple and $3$-connected.

Because $G_1\cap G_2=C$, we have $|E(G)|=|E(G_1)|+|E(G_2)|-3$.
Every triangle of $G$ is contained in $G_1$ or $G_2$, and $C$ is the unique triangle belonging to both subgraphs.
Thus $N(G,K_3)=N(G_1,K_3)+N(G_2,K_3)-1$, which yields $\delta(G)=\delta(G_1)+\delta(G_2)$.
Iterating this additivity gives $\delta(G)=\sum_A \delta(A)$.
Since an atom contains no separating triangles, every triangle in an atom $A$ is facial.
Therefore $|F(A)|=N(A,K_3)+q_A$.
By Euler's formula, $\delta(A)=v_A-4+q_A$.
Substituting it into $\delta(G)=\sum_A \delta(A)$ gives 
\(
\delta(G)=\sum_A(v_A-4+q_A).
\)

Define the atom adjacency graph $\mathcal A$ whose vertices are the atoms of $G$, with two atoms adjacent whenever they share a separating triangle from the decomposition.
We prove by induction on the number of separating triangles used in the decomposition that $\mathcal A$ is a tree.

If $G$ has no separating triangle, then $G$ itself is the only atom, so $\mathcal A$ is a single vertex.
Suppose the first decomposition step splits $G$ along a separating triangle $C$ into $G_1$ and $G_2$.
By the induction hypothesis, the atom graphs of $G_1$ and $G_2$ are trees.
Since $C$ is a facial triangle in both $G_1$ and $G_2$, it belongs to exactly one atom from each of these two trees.
We add an edge between these two atoms in $\mathcal A$.
Thus $\mathcal A$ is formed by connecting two trees with a single edge, so it is again a tree.

Now assume that every atom except $A_0$ is isomorphic to $K_4$.
Root the tree $\mathcal A$ at $A_0$.
We build $G$ incrementally starting from $A_0$, attaching the remaining atoms one at a time according to the tree order.
Each additional atom is a copy of $K_4$ that shares a triangular face with the graph constructed so far.
This attachment corresponds to adding one new vertex adjacent to the three vertices of that triangular face, which is precisely a stacking operation.
\end{proof}

\begin{lemma}
	\label{lem:small-delta}
	Let $G$ be a simple $3$-connected plane graph.
	Then
	\begin{enumerate}
		\item[(i)] $\delta(G)\ne1$;
		
		\item[(ii)] $\delta(G)=0$ if and only if $G$ is obtained from $K_4$ by a sequence of stackings;
		
		\item[(iii)] $\delta(G)=2$ if and only if {exactly one atom of $G$ is not isomorphic to $K_4$}, and this atom is $W_4$ or the octahedron.
	\end{enumerate}
\end{lemma}
\begin{proof}
	Let $A_1,\ldots,A_s$ be the atoms of $G$.
	By Lemma~\ref{lem:triangle-atoms},
	$\delta(G)=\sum_{i=1}^s\delta(A_i)$.
	For each $i$, let $q_i$ be the number of non-triangular faces of $A_i$.
	Then $\delta(A_i)=|V(A_i)|-4+q_i$.
	
(i)\ \ 
	If $|V(A_i)|\ge6$, then $\delta(A_i)\ge2$.
	If $|V(A_i)|=5$ and $q_i\ge1$, then $\delta(A_i)\ge2$.
	Suppose that $|V(A_i)|=5$ and $q_i=0$.
	Then $A_i$ is a triangulation with $9$ edges.
	By Lemma~\ref{lem:triangle-atoms}, every vertex of $A_i$ has degree
	at least three.
	Since the sum of the vertex degrees is $18$, there is a vertex of
	degree three.
	Its three neighbours form a separating triangle, a contradiction. If $|V(A_i)|=4$, then $A_i\cong K_4$ and $\delta(A_i)=0$.
	Therefore $\delta(G)=0$ or $\delta(G)\ge2$.	
	
(ii)\ \ 
	Suppose that $\delta(G)=0$.
	By Lemma~\ref{lem:triangle-atoms}, $\sum\limits_{A}(v_A-4+q_A)=0$.
	Thus every atom of $G$ is $K_4$. By the same argument as in the proof of
	Lemma~\ref{lem:triangle-atoms}, $G$ is obtained from $K_4$ by a
	sequence of stackings. Conversely, suppose that $G$ is obtained from $K_4$ by a sequence of
	stackings.
	Each stacking adds three edges and three triangles, so $\delta$ does not
	change.
	Since $\delta(K_4)=0$, we have $\delta(G)=0$.

(iii)\ \ Suppose that $\delta(G)=2$. Then exactly one atom is not $K_4$. Say
this atom is $A_1$. We have $\delta(A_1)=2$, and hence $|V(A_1)|-4+q_1=2$. Therefore,
\[
(|V(A_1)|,q_1)=(5,1)\quad\text{or}\quad(6,0).
\]
Suppose first that $|V(A_1)|=5$ and $q_1=1$. Then $A_1$ has exactly one non-triangular
face. By Euler's formula, $A_1$ has $|E(A_1)|-3$ faces. Thus it has $|E(A_1)|-4$ triangular
faces. Since the non-triangular face has length at least four,
\[
2|E(A_1)|\ge 3\bigl(|E(A_1)|-4\bigr)+4.
\]
Hence $|E(A_1)|\le 8$. Since $A_1$ is 3-connected, every vertex has degree at least three. Thus
$|E(A_1)|\ge 8$, and hence $|E(A_1)|=8$. The sum of the vertex degrees is 16. Thus one vertex
has degree four and the other four have degree three. The vertex of degree four is adjacent
to all other vertices. The other four vertices form a 4-cycle. Hence $A_1=W_4$.

Now suppose that $|V(A_1)|=6$ and $q_1=0$. Then $A_1$ is a triangulation on six vertices,
so it has twelve edges. No vertex has degree three. Otherwise, its three neighbours form
a separating triangle. Thus every vertex has degree at least four. Since the sum of the
vertex degrees is 24, every vertex has degree four. Hence each vertex is nonadjacent to
exactly one other vertex. It follows that $A_1$ is the octahedron. Conversely, if exactly one atom is $W_4$ or the octahedron and all
other atoms are $K_4$, then $\delta(G)=2$.
\end{proof}

Let $G$ be a simple $3$-connected plane graph, and let $G'$ be
obtained from $G$ by stacking a new vertex in a triangular face $f$.
This operation adds three edges and creates exactly three new
triangles, while preserving all existing triangles. Hence
$\delta(G')=\delta(G)$.

Let $x$ be the vertex of $G^*$ corresponding to $f$, with distinct
neighbours $a_1,a_2,a_3$. Then $(G')^*$ is obtained from $G^*$ by
replacing $x$ with a triangle $x_1x_2x_3$ and joining $x_i$ to $a_i$
for $i=1,2,3$. We call this operation a \emph{triangle expansion}
at $x$. Let $t_x$ denote the number of triangles of $G^*$ containing
$x$. The expansion removes these $t_x$ triangles and creates only
the triangle $x_1x_2x_3$: each $x_i$ has only one neighbour outside
the new triangle, and these neighbours are distinct.
The number of edges increases by three. Consequently,
\begin{equation}\label{eq:stacking-dual-deficit}
 \delta((G')^*)-\delta(G^*)
 =3-(1-t_x)=2+t_x\ge2.
\end{equation}
Thus stacking preserves the deficit of the primal graph and
increases that of its dual by at least two.

\begin{lemma}\label{lem:R-deficit}
Let $G$ be a simple $3$-connected plane graph.
\begin{enumerate}
 \item[(i)] If $G\not\cong K_4$, then $\Delta(G)\ge4$, with
 equality if and only if $G\cong W_4$.
 \item[(ii)] If $|E(G)|\ge12$, then $\Delta(G)\ge8$.
\end{enumerate}
\end{lemma}

\begin{proof}
Recall that the dual of a simple $3$-connected plane graph is
again simple and $3$-connected. Duality preserves both the number
of edges and $\Delta$. By Lemma~\ref{lem:small-delta}(i), the
deficits of such graphs are nonnegative and never equal one.

\smallskip
\noindent\textbf{(i)}
Suppose first that $\delta(G)=0$ or $\delta(G^*)=0$.
By duality, we may assume that $\delta(G)=0$.
Since $G\not\cong K_4$, Lemma~\ref{lem:small-delta}(ii) implies
that $G$ is obtained from $K_4$ by at least one stacking.
Let $G_1$ be the graph after the first stacking.
Since $K_4^*\cong K_4$ and every vertex of $K_4$ lies in three
triangles, \eqref{eq:stacking-dual-deficit} gives
\(
 \delta(G_1^*)=\delta(K_4^*)+2+3=5.
\)
Every subsequent stacking increases the dual deficit, so
$\Delta(G)=\delta(G^*)\ge5$.

Otherwise, both deficits are positive and hence at least two.
Thus $\Delta(G)\ge4$. If equality holds, then
$\delta(G)=\delta(G^*)=2$.
By Lemmas~\ref{lem:triangle-atoms}
and~\ref{lem:small-delta}(iii), $G$ is obtained by stackings from
either the octahedral graph or $W_4$.
In the former case, $G^*$ is obtained from the cube by triangle
expansions. Since the cube has $12$ edges and no triangles,
its deficit is $10$, giving $\delta(G^*)\ge10$, a contradiction.

Suppose that $G$ is obtained from $W_4$. If at least one stacking
is performed, then the first expands a degree-three vertex of
$W_4^*\cong W_4$. Every such vertex lies in exactly two triangles,
so the dual deficit after this stacking is
\(
 \delta(W_4^*)+2+2=6.
\)
Subsequent stackings only increase it, contradicting
$\delta(G^*)=2$. Hence $G\cong W_4$.
Conversely, $W_4$ is self-dual and $\delta(W_4)=2$, so
$\Delta(W_4)=4$.

\smallskip
\noindent\textbf{(ii)}
Suppose, to the contrary, that $|E(G)|\ge12$ and $\Delta(G)\le7$.
By duality, we may assume that $\delta(G)\le3$.
Lemma~\ref{lem:small-delta}(i) then gives
$\delta(G)\in\{0,2,3\}$.

\emph{Case 1: $\delta(G)=0$.}
By Lemma~\ref{lem:small-delta}(ii), $G$ is obtained from $K_4$
by stackings. Since each stacking adds three edges, at least two
are required. After the first stacking, the dual is the triangular
prism and has deficit $5$. Every vertex of the prism lies in
exactly one triangle, so the second stacking increases the dual
deficit by $3$. By \eqref{eq:stacking-dual-deficit}, further
stackings only increase it. Thus $\Delta(G)=\delta(G^*)\ge8$,
a contradiction.

\emph{Case 2: $\delta(G)=2$.}
By Lemmas~\ref{lem:triangle-atoms}
and~\ref{lem:small-delta}(iii), $G$ is obtained by stackings from
the octahedral graph or $W_4$. The octahedral case gives
$\delta(G^*)\ge10$ and is impossible.
Thus $G$ is obtained from $W_4$ by $r$ stackings, where
$8+3r=|E(G)|\ge12$ implies $r\ge2$.
The calculation in part~(i) gives dual deficit $6$ after the
first stacking, and each subsequent stacking increases it by at
least two. Hence $\delta(G^*)\ge8$, so $\Delta(G)\ge10$,
a contradiction.

\emph{Case 3: $\delta(G)=3$.}
By Lemmas~\ref{lem:triangle-atoms}
and~\ref{lem:small-delta}(i), the deficit is additive over the
atoms, and every atom other than $K_4$ has deficit at least two.
Consequently, $G$ has exactly one atom $A\not\cong K_4$, with
$\delta(A)=3$; all its other atoms are copies of $K_4$.
Rooting the atom decomposition at $A$ shows that $G$ is obtained
from $A$ by stackings.

Write $n_A=|V(A)|$ and $m_A=|E(A)|$, and let $q_A$ be the number
of nontriangular faces of $A$. Since $A$ is simple and
$3$-connected and is not $K_4$, we have $n_A\ge5$.
The identity $\delta(A)=n_A-4+q_A=3$ therefore gives
\(
 (n_A,q_A)\in\{(7,0),(6,1),(5,2)\}.
\)

If $(n_A,q_A)=(5,2)$, then $3$-connectivity gives $m_A\ge8$.
Euler's formula shows that $A$ has $m_A-3$ faces, of which
$m_A-5$ are triangular. Counting face lengths gives
\(
 2m_A\ge3(m_A-5)+2\cdot4=3m_A-7,
\)
so $m_A\le7$, a contradiction.

If $(n_A,q_A)=(7,0)$, then $A$ is a triangulation with $15$
edges, and $A^*$ is a simple $3$-connected cubic graph on
$10$ vertices. Its triangles are pairwise vertex-disjoint.
Indeed, cubicity forces any two distinct intersecting triangles
to share an edge; write them as $uvw$ and $uvz$.
The vertices $u$ and $v$ have no neighbours outside
$\{u,v,w,z\}$, so deleting $w$ and $z$ isolates the edge $uv$,
contradicting $3$-connectivity.
Thus $N(A^*,K_3)\le\lfloor10/3\rfloor=3$, and hence
\(
 \delta(A^*)=15-2-N(A^*,K_3)\ge10.
\)
Since triangle expansions increase the deficit,
$\delta(G^*)\ge\delta(A^*)\ge10$, again a contradiction.

It remains to consider $(n_A,q_A)=(6,1)$.
Let $\ell$ be the length of the unique nontriangular face of $A$.
Its boundary is a cycle, so $4\le\ell\le6$.
Euler's formula and the sum of the face lengths give
$m_A=15-\ell$.
If $\ell=6$, choose this face as the outer face.
Its boundary cycle $C$ contains all vertices of $A$, and
$m_A=9$ guarantees a chord $uv$ of $C$.
Both components of $C-\{u,v\}$ are nonempty.
All chords lie inside $C$, so an edge joining these components
would cross $uv$. Hence $A-\{u,v\}$ is disconnected,
contradicting $3$-connectivity.
Therefore $\ell\in\{4,5\}$ and $m_A\in\{10,11\}$.

We claim that $\delta(A^*)\ge3$.
Lemma~\ref{lem:small-delta}(i) excludes the value one.
If $\delta(A^*)=0$, Lemma~\ref{lem:small-delta}(ii) gives
$|E(A^*)|=6+3r$ for some integer $r\ge0$, contrary to
$|E(A^*)|=m_A\in\{10,11\}$.
If $\delta(A^*)=2$, Lemma~\ref{lem:small-delta}(iii) gives
$|E(A^*)|=8+3r$ or $12+3r$ for some integer $r\ge0$.
Thus $|E(A^*)|=11$, and $A^*$ is obtained from $W_4$ by
exactly one stacking, so it has six vertices. However, Euler's
formula gives $|V(A^*)|=m_A-n_A+2=7$, a contradiction.
This proves the claim.

Finally, $m_A<12\le|E(G)|$ implies that at least one stacking is
performed in obtaining $G$ from $A$. Hence, 
\(
 \delta(G^*)\ge\delta(A^*)+2\ge5,
\) and so \( 
 \Delta(G)\ge3+5=8,
\)
a contradiction.
\end{proof}

\subsection{The sharp bound and extremal quadrangulations}
\label{sec:main-proof}
We now prove that $\mu(Q)\ge6$ for every $3$-connected plane
quadrangulation $Q$ with at least $14$ vertices and characterize
the cases of equality.

\begin{lemma}\label{lem:K4-loss}
Let $H$ be a simple $2$-connected plane graph whose dual $H^*$
is simple, and let $\mathcal T$ be its SPQR tree.
Let $\rho$ be an $R$-node of $\mathcal T$ with $J_\rho\cong K_4$.
If $c_\rho\ge1$, then $\ell_\rho\ge2$.
If $c_\rho\ge2$, then $\ell_\rho\ge3$.
\end{lemma}

\begin{proof}
Suppose that $c_\rho\ge1$, and let $e$ be a virtual edge of
$J_\rho$ corresponding to a core edge $\rho x$.
By definition, $e$ is $H$-active only if $x$ is a $P$-node in $U$,
whereas its dual edge $e^*$ is $H^*$-active only if $x$ is an
$S$-node in $U$.
Thus at least one of $e$ and $e^*$ is inactive.
Since $J_\rho\cong J_\rho^*\cong K_4$ and every edge of $K_4$
lies in exactly two triangles, $\ell_\rho\ge2$.

Now suppose that $c_\rho\ge2$.
Choose distinct virtual edges $e_1,e_2$ of $J_\rho$ and an
inactive edge from each pair $\{e_i,e_i^*\}$.
If the chosen edges lie in different skeletons, they account
for four triangles counted by $\ell_\rho$.
If they lie in the same skeleton, they are distinct edges of
$K_4$ and together belong to at least three triangles.
In either case, $\ell_\rho\ge3$.
\end{proof}

\begin{proposition}\label{prop:core}
Let $Q$ be a $3$-connected plane quadrangulation with $n\ge14$
vertices and bipartition $V(Q)=X\cup Y$.
Set $H:=Q_X$ and $m:=|E(H)|$, and let $\mathcal T$ be the
SPQR tree of $H$.
Then $\mu(Q)\ge6$, with equality if and only if the core of
$\mathcal T$ is a path $R_0X_1\cdots X_kR_{k+1}$, where
$k\ge2$, satisfying the following conditions:
\begin{itemize}
\item the skeletons of $R_0$ and $R_{k+1}$ are both $K_4$,
each with five real edges and one virtual edge;
\item each $X_i$ belongs to $U$, with $d_{X_i}=3$ and $c_{X_i}=2$;
\item $X_1,\ldots,X_k$ alternate between $S$-nodes and $P$-nodes.
\end{itemize}
In this case, $m=10+k$ and $n=12+k$.
\end{proposition}

\begin{proof}
By Lemma~\ref{lem:radial}, $H^*=Q_Y$ and $m=n-2$.
We use the identity \eqref{eq:potential}, whose
summands indexed by $U\cup W$ are nonnegative.

Suppose first that the core consists of a single node $\rho$.
By Lemma~\ref{lem:core-leaves}, $\rho$ is an $R$-node.
Its skeleton is $H$, and all edges of $J_\rho$ and $J_\rho^*$
are real and hence active.
Thus $\ell_\rho=0$.
Since $m=n-2\ge12$, Lemma~\ref{lem:R-deficit}(ii) gives
\(
\mu(Q)=\Delta(J_\rho)\ge8.
\)

Henceforth suppose that the core has at least two nodes.
Every leaf is an $R$-node by Lemma~\ref{lem:core-leaves}, so
$|\mathcal R|\ge2$.
Every $R$-node has positive core degree.
If $J_\rho\cong K_4$, then
$\Delta(J_\rho)+\ell_\rho=\ell_\rho\ge2$ by
Lemma~\ref{lem:K4-loss}.

We claim that $\Delta(J_\rho)+\ell_\rho\ge5$ whenever
$J_\rho\not\cong K_4$.
If $J_\rho\not\cong W_4$, this follows from
Lemma~\ref{lem:R-deficit}(i) and the integrality of $\Delta(J_\rho)$.
If $J_\rho\cong W_4$, choose a virtual edge $e$ of $J_\rho$.
As in the proof of Lemma~\ref{lem:K4-loss}, at least one of
$e$ and $e^*$ is inactive.
Since $W_4^*\cong W_4$ and every edge of $W_4$ lies in a
triangle, $\ell_\rho\ge1$.
Together with $\Delta(W_4)=4$, this proves the claim.
Consequently, if some $R$-skeleton is not $K_4$, its contribution
to \eqref{eq:potential} is at least five, and that of another
$R$-node is at least two. Hence $\mu(Q)\ge7$.

We may therefore assume that every $R$-skeleton is $K_4$.
If $|\mathcal R|\ge4$, then $\mu(Q)\ge8$.
Suppose that $|\mathcal R|=3$.
If some $R$-node has core degree at least two,
Lemma~\ref{lem:K4-loss} gives $\mu(Q)\ge3+2+2=7$.
Otherwise, all three $R$-nodes are leaves.
The core then has exactly three leaves and hence a node $x$
with $c_x\ge3$.
This node belongs to $U\cup W$.
If $x\in U$, then $c_x-1-z_x\ge1$; if $x\in W$, then
$2c_x-4-z_x\ge1$.
Thus \eqref{eq:potential} again gives $\mu(Q)\ge7$.

It remains to consider the case in which the core has exactly
two $R$-nodes, both with skeleton $K_4$.
Since every leaf is an $R$-node, the core is a path
\(
R_0X_1\cdots X_kR_{k+1}.
\) 
Each end skeleton has one virtual edge and five real edges.
Every internal node $X_i$ is an $S$- or $P$-node with
$c_{X_i}=2$.
As every node in $W$ has core degree at least three,
$X_i\in U$ and $d_{X_i}=3$.
Thus its skeleton has one real edge and two virtual edges.
Consecutive internal nodes have different types, since no two
$S$-nodes and no two $P$-nodes are adjacent.
Counting real edges gives
\(
m=10+k, n=m+2=12+k.
\) 
In particular, $k\ge2$.

For either end node $\rho$, let $e$ be its unique virtual edge.
Its neighbour in the core belongs to $U$, so exactly one of
$e$ and $e^*$ is active.
All other edges of $J_\rho$ and $J_\rho^*$ are real and hence
active. Therefore
\(
\ell_{R_0}=\ell_{R_{k+1}}=2.
\) 
The edge corresponding to $R_0X_1$ is inactive in $J_{X_1}$
if $X_1$ is an $S$-node, and its dual is inactive in
$J_{X_1}^*$ if $X_1$ is a $P$-node.
Lemma~\ref{lem:triangle-home} therefore gives $z_{X_1}=0$.
Similarly, $z_{X_k}=0$.
For $2\le i\le k-1$, both neighbours of $X_i$ belong to $U$
and have the opposite type.
Thus all three edges of $J_{X_i}$ are $H$-active if $X_i$ is
an $S$-node, and all three edges of $J_{X_i}^*$ are
$H^*$-active if $X_i$ is a $P$-node.
By Lemma~\ref{lem:triangle-home}, $z_{X_i}=1$.
Since $\Delta(K_4)=0$, \eqref{eq:potential} gives
\(
\mu(Q)
=\ell_{R_0}+\ell_{R_{k+1}}+\sum_{i=1}^{k}(1-z_{X_i})
=2+2+1+1=6.
\)
This calculation applies to every core satisfying the stated
conditions, whereas all other cases give $\mu(Q)\ge7$.
The bound and the equality characterization follow.
\end{proof}

Proposition~\ref{prop:core} gives
$s_3(Q)=2n-8-\mu(Q)\le2n-14$ for $n\ge14$.
We next translate the equality conditions into the construction
of $\mathcal G^*$.
For an edge $e=uv$ of a $2$-connected plane graph $G$, let $f,g$
be its incident faces. The corresponding face of
$\operatorname{Rad}(G)$ is bounded by the $4$-cycle
$C(e):=ufvgu$.

\begin{lemma}\label{lem:radial-gluing}
Let $G_1,G_2$ be disjoint $2$-connected plane graphs, and let
$e_1=u_1v_1\in E(G_1)$ and $e_2=u_2v_2\in E(G_2)$.
Let $f_1,f_2$ be the faces incident with $e_1$, and let
$g_1,g_2$ be the faces incident with $e_2$.
Form a plane graph $M$ by identifying $u_1$ with $u_2$ and
$v_1$ with $v_2$, and deleting $e_1,e_2$, with the embeddings
chosen so that $f_i$ and $g_i$ merge into one face for $i=1,2$.
Then $\operatorname{Rad}(M)$ is obtained by gluing
$\operatorname{Rad}(G_1)$ and $\operatorname{Rad}(G_2)$ along
the facial cycles
\[
C(e_1)=u_1f_1v_1f_2u_1,\qquad
C(e_2)=u_2g_1v_2g_2u_2,
\]
identifying $u_1$ with $u_2$, $v_1$ with $v_2$, and $f_i$
with $g_i$ for $i=1,2$, together with the corresponding
boundary edges.
\end{lemma}

\begin{proof}
For $i=1,2$, removing $e_1$ and $e_2$ from the boundaries
of $f_i$ and $g_i$, respectively, leaves two paths.
After the end identifications, these paths form the
boundary of a face $h_i$ of $M$.
All other faces remain unchanged.

Thus the vertices of $\operatorname{Rad}(M)$ are obtained by
identifying $u_1$ with $u_2$, $v_1$ with $v_2$, and each
pair $f_i,g_i$ with the face vertex $h_i$.
The radial edges $u_1f_i$ and $u_2g_i$ become the single
edge $u_1h_i$, and similarly at $v_1,v_2$.
All remaining vertex--face incidences are preserved. These are precisely the identifications of the corresponding
vertices and edges of $C(e_1)$ and $C(e_2)$.
\end{proof}

\begin{lemma}\label{lem:expansions}
Let $G$ be a $2$-connected plane graph, and let $e=uv$ be
incident with the faces $f,g$.
\begin{enumerate}
\item[(i)] Let $G_1$ be obtained from $G$ by subdividing $e$
with a new vertex $z$.
Then $\operatorname{Rad}(G_1)$ is obtained from
$\operatorname{Rad}(G)$ by adding $z$ inside the face bounded
by $C(e)$ and joining it to $f$ and $g$.
\item[(ii)] Let $G_2$ be obtained from $G$ by replacing $e$
with two parallel edges bounding a new face $h$.
Then $\operatorname{Rad}(G_2)$ is obtained from
$\operatorname{Rad}(G)$ by adding $h$ inside the face bounded
by $C(e)$ and joining it to $u$ and $v$.
\end{enumerate}
Both $G_1$ and $G_2$ are $2$-connected.
\end{lemma}

\begin{proof}
In (i), the only new vertex--face incidences are those of $z$
with $f$ and $g$.
Hence the radial graph gains exactly the vertex $z$ and the
edges $zf,zg$ inside $C(e)$.
In (ii), the only new face is $h$, whose boundary consists
of the two parallel edges with ends $u,v$.
Hence the radial graph gains exactly the vertex $h$ and the
edges $hu,hv$ inside $C(e)$.

Adding a parallel edge preserves $2$-connectivity, so $G_2$
is $2$-connected.
Since $e$ is not a bridge, $G_1-z=G-e$ is connected.
For each $x\in V(G)$, the graph $G_1-x$ is obtained from
the connected graph $G-x$ either by subdividing an edge
or by attaching $z$ to the remaining end of $e$.
Thus $G_1$ is also $2$-connected.
\end{proof}

We call the operations on $G$ in
Lemma~\ref{lem:expansions}(i) and (ii) an \emph{$S$-expansion}
and a \emph{$P$-expansion}, respectively.

\begin{proposition}\label{prop:radial-characterization}
Let $H$ be a $2$-connected plane graph whose SPQR tree has
core $R_0X_1\cdots X_kR_{k+1}$, where $k\ge1$.
Suppose that the end skeletons are both $K_4$, each internal
node $X_i$ has one real edge and two virtual edges, and
$X_1,\ldots,X_k$ alternate between $S$-nodes and $P$-nodes.
Then $\operatorname{Rad}(H)\cong G_k$ for some graph $G_k$
in Construction~\ref{con}.
Conversely, for any graph $G_k$ in Construction~\ref{con}
with bipartition $V(G_k)=X\cup Y$, the graph $H:=(G_k)_X$
is $2$-connected and its SPQR tree has the core and skeletons
described above.
\end{proposition}

\begin{proof}
Let $J_0$ be the skeleton at $R_0$.
Thus $J_0\cong K_4$, with five real edges and one virtual
edge $e_0$.
For $1\le i\le k$, glue the skeleton at $X_i$ to $J_{i-1}$
along $e_{i-1}$ and its twin, identifying corresponding
ends and deleting both edges.
Denote the resulting graph by $J_i$ and its unique virtual
edge by $e_i$.
If $X_i$ is an $S$-node, this subdivides $e_{i-1}$;
if $X_i$ is a $P$-node, it replaces $e_{i-1}$ with two
parallel edges.
In either case, one new edge is real and the other is $e_i$.
By Lemma~\ref{lem:expansions}, every $J_i$ is $2$-connected.

For $0\le i\le k$, designate the face of
$\operatorname{Rad}(J_i)$ bounded by $C(e_i)$ as active.
Each step inserts a vertex $w_i$ into the active face of
$\operatorname{Rad}(J_{i-1})$ and joins it to two opposite
vertices of that face.
Since $\operatorname{Rad}(K_4)$ is $K_{4,4}$ minus a perfect
matching, it is the cube.
Label $C(e_0)=v_0v_1v_2v_3v_0$ so that $w_1$ is joined
to $v_1,v_3$.

An $S$-expansion adds a primal vertex to the radial graph,
whereas a $P$-expansion adds a face vertex.
Thus $w_1,\ldots,w_k$ alternate between the two partite sets.
Each $w_i$ lies on the boundary of the new active face.
Consequently, for $i\ge2$, the vertex $w_i$ is joined to
$w_{i-1}$ and to the vertex opposite $w_{i-1}$ on the current
active face.
These are exactly the insertion steps in Construction~\ref{con}.
Finally, gluing the skeleton at $R_{k+1}$ to $J_k$ along
$e_k$ and its twin corresponds, by Lemma~\ref{lem:radial-gluing},
to attaching a second cube along $C(e_k)$.
Hence $\operatorname{Rad}(H)\cong G_k$ for a graph $G_k$
in Construction~\ref{con}.

Conversely, let $G_k$ be a graph in Construction~\ref{con},
with bipartition $V(G_k)=X\cup Y$, and set $H:=(G_k)_X$.
Since $G_k$ is a $3$-connected plane quadrangulation,
Lemma~\ref{lem:radial} implies that $H$ is simple and
$2$-connected.
Let $Q_0$ be the left cube, with its attachment face
designated as active.
For $1\le i\le k$, obtain $Q_i$ by inserting $w_i$ into
the active face of $Q_{i-1}$ as in Construction~\ref{con},
and designate the prescribed new face as active.

For $0\le i\le k$, let $J_i$ have vertex set
$X\cap V(Q_i)$ and one edge for the $X$-diagonal of each
face of $Q_i$, retaining parallel edges when they occur.
Let $e_i$ be the diagonal of the active face, and mark
$e_i$ as virtual and all other edges of $J_i$ as real.
Every nonactive face of $Q_i$ remains a face of $G_k$,
so all these real edges belong to $H$.
We have $J_0\cong K_4$.
If $w_i\in X$, the two new diagonals replace $e_{i-1}$
by a path of length two through $w_i$; if $w_i\in Y$,
they replace $e_{i-1}$ by two parallel edges.
These are an $S$-expansion and a $P$-expansion, respectively.
In either case, one new edge is real and the other is $e_i$.
Since $w_1,\ldots,w_k$ alternate between $X$ and $Y$,
the expansions alternate between the two types.
Every $J_i$ is $2$-connected by Lemma~\ref{lem:expansions}.

The $X$-diagonal graph of the right cube is another copy
of $K_4$.
Thus attaching that cube to the active face of $Q_k$
corresponds to identifying the ends of $e_k$ with those
of an edge $e_k'$ of this copy of $K_4$, and deleting
$e_k,e_k'$.
The resulting graph is $H$.

We now reverse this construction by successive splits.
First split $H$ into $J_k$ and the right copy of $K_4$,
restoring $e_k,e_k'$ as twin virtual edges.
Next, for $i=k,k-1,\ldots,1$, partition $E(J_i)$ into
$E(J_{i-1})\setminus\{e_{i-1}\}$ and the two edges
introduced at step $i$.
The corresponding subgraphs meet exactly at the ends
of $e_{i-1}$, and each has at least two edges.
Splitting according to this partition yields $J_{i-1}$
and a skeleton $N_i$.
For an $S$-expansion, $N_i$ is a triangle;
for a $P$-expansion, it consists of three parallel edges.
In either case, $N_i$ has one real edge and two virtual
edges: $e_i$ and the new virtual edge whose twin is
$e_{i-1}$.
All existing twin relations are retained.

After these splits, the remaining graph is $J_0\cong K_4$.
The resulting skeletons form the path
$R_0X_1\cdots X_kR_{k+1}$, with $N_i$ at $X_i$ and
a copy of $K_4$ at each end.
Since the $S$- and $P$-nodes alternate, no further merging
is required.
Adding a $Q$-node for each real edge therefore gives an
SPQR tree of $H$.
Lemma~\ref{lem:spqr-uniqueness} completes the proof.
\end{proof}

\begin{corollary}\label{cor:Gk-count}
Let $G_k$ be a graph in Construction~\ref{con}, and let
$n=|V(G_k)|=k+12$. Then $s_3(G_k)=2n-14$ for $k\ge2$, and $s_3(G_1)=13$.
\end{corollary}

\begin{proof}
Let $V(G_k)=X\cup Y$ be its bipartition and set $H:=(G_k)_X$.
By Proposition~\ref{prop:radial-characterization}, the core
of the SPQR tree of $H$ is $R_0X_1\cdots X_kR_{k+1}$
with the stated skeletons.
For $k\ge2$, Proposition~\ref{prop:core} gives $\mu(G_k)=6$,
and hence $s_3(G_k)=2n-14$.

Suppose that $k=1$, so the core is $R_0X_1R_2$.
For each $\rho\in\{R_0,R_2\}$, the skeleton $J_\rho$ is
$K_4$ and has a unique virtual edge $e$.
Since $X_1\in U$ is an $S$- or $P$-node, exactly one of
$e,e^*$ is inactive; all other edges of $J_\rho$ and
$J_\rho^*$ are real and active.
Thus $\Delta(J_\rho)=0$ and $\ell_\rho=2$.
Both core neighbours of $X_1$ are $R$-nodes, so its two
virtual edges are inactive on the relevant side.
By Lemma~\ref{lem:triangle-home}, $z_{X_1}=0$.
Since $c_{X_1}=2$, \eqref{eq:potential} gives
\(
\mu(G_1)=2+2+(2-1-0)=5.
\)
Therefore $s_3(G_1)=2\cdot13-8-5=13$.
\end{proof}

We next determine the maximum number of triangles in optimal
$1$-planar graphs of orders $10\le n\le13$.

\begin{proposition}\label{prop:small-triangles}
Let $G$ be an optimal $1$-planar graph with $n$ vertices,
where $10\le n\le13$. Then
\[
N(G,K_3)\le
\begin{cases}
40, & n=10,\\
45, & n=11,\\
50, & n=12,\\
57, & n=13,
\end{cases}
\]
with equality if and only if $G\cong\widehat Q$, where $Q$ is
$\operatorname{Rad}(W_4)$,
$\operatorname{Rad}(C_3\square K_2)$,
$\operatorname{Rad}(W_5)$, or a graph $G_1$ from
Construction~\ref{con}, respectively.
\end{proposition}

\begin{proof}
By Proposition~\ref{prop:vertex-face},
$G\cong\widehat{\operatorname{Rad}(H)}$ for a simple
$2$-connected plane graph $H$ with simple dual.
Set $Q:=\operatorname{Rad}(H)$, and let $\mathcal T^\circ$
be the core of the SPQR tree of $H$, with set $\mathcal R$ of
$R$-nodes.
By Corollary~\ref{cor:triangles} and the definition of $\mu(Q)$,
\(
N(G,K_3)=6n-16-\mu(Q).
\)
Lemma~\ref{lem:core-leaves} and \eqref{eq:spqr-size} imply
$1\le|\mathcal R|\le2$.

Suppose first that $|\mathcal R|=2$.
By Lemma~\ref{lem:core-leaves}, the core is a path whose
ends are the two $R$-nodes.
Every $R$-skeleton other than $K_4$ has at least eight edges,
which would give $n\ge14$ by \eqref{eq:spqr-size}.
Thus both end skeletons are $K_4$.
Every internal node has core degree two, so it belongs to
$U$ and has one real edge and two virtual edges.
By \eqref{eq:spqr-size}, the path has $n-12$ internal
nodes. Hence $n\in\{12,13\}$.

If $n=12$, the two $R$-nodes are adjacent.
For each such node $\rho$, its unique virtual edge and
the corresponding dual edge are both inactive, while all
other edges are real and active.
Since $J_\rho\cong J_\rho^*\cong K_4$, we obtain
$\ell_\rho=4$. By \eqref{eq:potential}, $\mu(Q)=8$, so $N(G,K_3)=48$.
If $n=13$, the core has one internal $S$- or $P$-node.
Proposition~\ref{prop:radial-characterization} gives
$Q\cong G_1$, and Corollary~\ref{cor:Gk-count} gives
$s_3(Q)=13$.
Thus $N(G,K_3)=57$.

Now suppose that $|\mathcal R|=1$.
By Lemma~\ref{lem:core-leaves}, the core consists of a
single $R$-node.
Thus $H$ is a simple $3$-connected plane graph with $n-2$
edges, and all edges of its skeleton and dual skeleton
are active.
By \eqref{eq:potential}, $\mu(Q)=\Delta(H)$, so
\(
N(G,K_3)=6n-16-\Delta(H).
\)
Since $H$ and $H^*$ are simple and each has at least eight
edges, both have at least five vertices.

If $n=10$, Euler's formula gives
$|V(H)|=|V(H^*)|=5$.
The degree sequence of $H$ is $(4,3,3,3,3)$, so
$H\cong W_4$.
Since $W_4^*\cong W_4$ and $N(W_4,K_3)=4$, we have
$\Delta(H)=4$ and $N(G,K_3)=40$.

If $n=11$, then $|E(H)|=9$, and Euler's formula gives
\(
\{|V(H)|,|V(H^*)|\}=\{5,6\}.
\)
Since $\operatorname{Rad}(H)\cong\operatorname{Rad}(H^*)$,
we may interchange $H$ and $H^*$ and assume that $|V(H^*)|=5$. Then $H^*\cong K_5-e$ and $H\cong C_3\square K_2$.
These graphs contain seven and two triangles, respectively. Thus $\Delta(H)=2\cdot9-4-7-2=5$ and $N(G,K_3)=45$.

Suppose that $n=12$, so $|E(H)|=10$.
A simple plane graph on five vertices has at most nine
edges. Euler's formula therefore gives
$|V(H)|=|V(H^*)|=6$.
The degree sequence of either graph is
$(5,3,3,3,3,3)$ or $(4,4,3,3,3,3)$.
If either graph has a vertex of degree five, deleting that
vertex leaves a $2$-regular simple graph on five vertices,
which must be $C_5$.
Thus that graph is $W_5$; since $W_5^*\cong W_5$,
we obtain $H\cong W_5$.
In this case, $\Delta(H)=6$ and $N(G,K_3)=50$.

Otherwise, both $H$ and $H^*$ have degree sequence
$(4,4,3,3,3,3)$.
Neither graph has a nonfacial triangle.
Indeed, such a triangle would be separating, with one
of the remaining three vertices on one side and two on
the other.
Decomposing along it would give two simple $3$-connected
plane graphs on four and five vertices, by
Lemma~\ref{lem:triangle-atoms}.
Their edge counts would imply
\(
10\ge6+8-3=11,
\)
a contradiction.
By duality, each of $H$ and $H^*$ has exactly four
triangular faces.
Since every triangle is facial,
$N(H,K_3)=N(H^*,K_3)=4$.
Hence $\Delta(H)=8$ and $N(G,K_3)=48$.

Finally, suppose that $n=13$, so $|E(H)|=11$.
The deficits $\delta(H),\delta(H^*)$ are nonnegative by
Lemma~\ref{lem:triangle-atoms}.
Lemma~\ref{lem:small-delta}(i),(ii) excludes the values
one and zero, respectively, since a graph of deficit zero
has $6+3s$ edges for some integer $s\ge0$.
Thus $\delta(H),\delta(H^*)\ge2$.
If either deficit equals two, interchange $H$ and $H^*$,
if necessary, so that $\delta(H)=2$.
Lemmas~\ref{lem:triangle-atoms} and~\ref{lem:small-delta}(iii),
together with $|E(H)|=11$, imply that $H$ is obtained
from $W_4$ by one stacking.
Since $W_4^*\cong W_4$ and every vertex of $W_4$ lies
in at least two triangles, \eqref{eq:stacking-dual-deficit}
gives $\delta(H^*)\ge6$.
Thus $\Delta(H)\ge8$.
If both deficits are at least three, then $\Delta(H)\ge6$.
In either case,
\(
N(G,K_3)\le6n-22=56<57.
\)
The stated bounds and equality cases follow.
\end{proof}

\begin{proof}[\bf Proof of Theorem~\ref{thm}]
Let $Q$ be a $3$-connected plane quadrangulation with
$n\ge14$ vertices and bipartition $V(Q)=X\cup Y$.
By Proposition~\ref{prop:core},
\(
s_3(Q)=2n-8-\mu(Q)\le2n-14.
\)
If equality holds, then $\mu(Q)=6$.
Proposition~\ref{prop:core} describes the core of the SPQR
tree of $Q_X$ as a path with two end skeletons $K_4$ and
$k=n-12\ge2$ internal nodes, each with one real edge and
two virtual edges, alternating between $S$-nodes and $P$-nodes.
Since $Q\cong\operatorname{Rad}(Q_X)$ by
Lemma~\ref{lem:radial}, Proposition~\ref{prop:radial-characterization}
implies that $Q\cong G_k$ for a graph $G_k$ in
Construction~\ref{con}.
Conversely, every such $G_k$ is a $3$-connected plane
quadrangulation, and Corollary~\ref{cor:Gk-count} gives
$s_3(G_k)=2n-14$.
\end{proof}

Let $G$ be an optimal $1$-planar graph with $n\ge14$ vertices.
By Proposition~\ref{prop:vertex-face},
$G\cong\widehat Q$, where $Q=\operatorname{Rad}(H)$ for
a simple $2$-connected plane graph $H$ with simple dual.
Theorem~\ref{thm} and Corollary~\ref{cor:triangles} give
\(
N(G,K_3)=4n-8+s_3(Q)\le6n-22,
\)
with equality if and only if $Q$ is isomorphic to a graph
in $\mathcal G^*$.
Together with Proposition~\ref{prop:small-triangles}, this
proves the assertions of Theorem~\ref{thm:cliques}
concerning $K_3$.

\section{Proof of Theorem~\ref{thm:cliques}: Counting copies of $K_4$ and $K_5$}
\label{sec:larger-cliques}

For a plane graph $G$, let $v_3(G)$ denote the number of vertices
of degree three, and let $f_3(G)$ denote the number of triangular
faces.

\begin{lemma}\label{lem:R-full-cliques}
Let $G$ be a simple $3$-connected plane graph.
Then $\widehat{\operatorname{Rad}(G)}$ is $K_5$-free, and
\begin{equation}\label{eq:R-full-four}
N(\widehat{\operatorname{Rad}(G)},K_4)
=|E(G)|+v_3(G)+f_3(G)+N(G,K_4)+N(G^*,K_4).
\end{equation}
If $G\cong K_4$, then $N(\widehat{\operatorname{Rad}(G)},K_4)=16$.
If $G\not\cong K_4$, then
\(
N(\widehat{\operatorname{Rad}(G)},K_4)\le3|E(G)|-8.
\)
Furthermore, if $|E(G)|\ge10$, then
\(
N(\widehat{\operatorname{Rad}(G)},K_4)\le3|E(G)|-10.
\)
\end{lemma}

\begin{proof}
We first prove \eqref{eq:R-full-four} by classifying the copies
of $K_4$ according to their numbers of primal and face vertices.
For each edge $uv\in E(G)$ with incident faces $f,g$,
the set $\{u,v,f,g\}$ induces a $K_4$ in
$\widehat{\operatorname{Rad}(G)}$.
Since $G-\{u,v\}$ is connected,
Lemma~\ref{lem:common-faces} implies that $f,g$ are the only
faces whose boundaries contain both $u$ and $v$.
Thus exactly $|E(G)|$ copies have two primal and two face vertices.

Suppose that $\{u,v,w,f\}$ induces a $K_4$, where
$u,v,w\in V(G)$ and $f\in F(G)$.
Then $uvw$ is a triangle and all three vertices lie on the
boundary of $f$.
By Lemma~\ref{lem:common-faces}, all three edges of this triangle
lie on that boundary.
Since every face boundary is a cycle, $f$ is bounded by $uvw$.
Hence exactly $f_3(G)$ copies have three primal vertices and
one face vertex.
By duality, exactly $v_3(G)$ copies have one primal vertex and
three face vertices.
The remaining copies lie entirely in $G$ or $G^*$, proving
\eqref{eq:R-full-four}.

We next show that $\widehat{\operatorname{Rad}(G)}$ is $K_5$-free.
Neither $G$ nor $G^*$ contains a $K_5$, so, by duality, it
suffices to consider a $K_5$ with three or four primal vertices.
In the latter case, all four vertices lie on the boundary of
one face, whereas any three of them already force that boundary
to be a triangle by the preceding argument.
In the former case, two distinct faces have the same triangular
boundary, forcing $G=C_3$.
Both cases are impossible.

If $G\cong K_4$, then $G^*\cong K_4$, and
\eqref{eq:R-full-four} gives the value $16$.
Now suppose that $G\not\cong K_4$, then $|V(G)|\ge5$.
Every copy of $K_4$ in $G$ contains at most one vertex of degree
three in $G$.
Indeed, if two vertices of such a copy both had degree three,
deleting the other two vertices would separate them from the
vertices outside the copy, contrary to $3$-connectivity.

We claim that
\[
N(G,K_4)+v_3(G)\le2|V(G)|-6,
\]
with equality only if $|E(G)|\le9$.
Suppose first that at most two vertices have degree at least
four. Then $N(G,K_4)=0$.
If $|V(G)|\ge6$, then $v_3(G)\le |V(G)|\le2|V(G)|-6$.
If $|V(G)|=5$, the parity of the degree sum gives $v_3(G)\le4=2|V(G)|-6$.
Equality is possible only when
$(|V(G)|,v_3(G))=(6,6)$ or $(5,4)$, giving $|E(G)|=9$ or $8$,
respectively.

Suppose instead that at least three vertices have degree at
least four.
Then $v_3(G)\le |V(G)|-3$, and Lemma~\ref{lem:planar-K4} gives
$N(G,K_4)\le |V(G)|-3$.
If equality holds in the claimed bound, exactly three vertices
have degree at least four and
\(
N(G,K_4)=v_3(G)=|V(G)|-3.
\)
Every copy of $K_4$ contains these three vertices and one
vertex of degree three.
Since a vertex of degree three belongs to at most one copy
of $K_4$, every such vertex is adjacent to all three vertices
of degree at least four.
There are at most two vertices of degree three, since otherwise
$G$ contains $K_{3,3}$.
As $|V(G)|\ge5$, this forces $|V(G)|=5$ and $|E(G)|=9$.
The claim follows.

Consequently, if $|E(G)|\ge10$, integrality gives
\(
N(G,K_4)+v_3(G)\le2|V(G)|-7.
\)
The same bounds apply to $G^*$, which is simple,
$3$-connected, and not isomorphic to $K_4$.
Using $f_3(G)=v_3(G^*)$ and Euler's formula in
\eqref{eq:R-full-four}, we obtain
\[
N(\widehat{\operatorname{Rad}(G)},K_4)
\le |E(G)|+2|V(G)|+2|V(G^*)|-12
=3|E(G)|-8.
\]
If $|E(G)|\ge10$, applying the stronger bound to both $G$
and $G^*$ gives $N(\widehat{\operatorname{Rad}(G)},K_4)
\le3|E(G)|-10$.
\end{proof}

A \emph{tree decomposition} of a graph $G$ is a pair
$(T,(B_x)_{x\in V(T)})$, where $T$ is a tree and each
$B_x\subseteq V(G)$ is a \emph{bag}, satisfying:
\begin{enumerate}
\item[(i)] every edge of $G$ has both ends in some bag;
\item[(ii)] for every $v\in V(G)$, the nodes $x$ with
$v\in B_x$ induce a nonempty subtree of $T$.
\end{enumerate}

For the remainder of this section, let $G$ be an optimal
$1$-planar graph with $n$ vertices.
By Proposition~\ref{prop:vertex-face}, we may assume that
$G=\widehat{\operatorname{Rad}(H)}$, where $H$ is a simple
$2$-connected plane graph with simple dual.
Let $\mathcal T$ be the SPQR tree of $H$, with core
$\mathcal T^\circ$.
For each $x\in V(\mathcal T^\circ)$, let $J_x$ be its skeleton
with the embedding inherited from $H$.

During reconstruction, faces incident with twin virtual edges
are merged in pairs.
Each face of $J_x$ corresponds to a unique face of $H$, and
distinct faces of $J_x$ correspond to distinct faces of $H$.
For each $f\in F(H)$, let $y_f$ denote the corresponding vertex
of $G$.
Let $F_x$ be the set of vertices $y_f$ for which $f$ corresponds
to a face of $J_x$, and set
\(
B_x:=V(J_x)\cup F_x.
\)

\begin{lemma}\label{lem:vertex-face-bags}
The following statements hold.
\begin{enumerate}
\item[(i)]
$(\mathcal T^\circ,(B_x)_{x\in V(\mathcal T^\circ)})$
is a tree decomposition of $G$.
\item[(ii)]
Let $xy\in E(\mathcal T^\circ)$, and let $e_x,e_y$ be its
twin virtual edges in $J_x,J_y$, with common ends $u,v$.
Let $f,g$ be the faces of $H$ corresponding to the two faces
incident with $e_x$ in $J_x$.
Then
\[
B_x\cap B_y=\{u,v,y_f,y_g\},
\qquad
uy_f,uy_g,vy_f,vy_g\in E(G).
\]
Moreover, $uv\in E(G)$ if and only if $x$ or $y$ is a
$P$-node in $U$, and $y_fy_g\in E(G)$ if and only if $x$
or $y$ is an $S$-node in $U$.
\item[(iii)]
For every $R$-node $\rho$, the graph $G[B_\rho]$ is
isomorphic to a graph obtained from $\operatorname{Rad}(J_\rho)$
by adding diagonals to its faces.
The face corresponding to a real edge of $J_\rho$ receives
both diagonals, whereas the face corresponding to a virtual
edge receives at most one.
\end{enumerate}
\end{lemma}

\begin{proof}
\textup{(i)}
Let $uv\in E(H)$ have incident faces $f,g$.
The edge $uv$ is real in some skeleton $J_x$, whose incident
faces correspond to $f,g$.
Thus $\{u,v,y_f,y_g\}\subseteq B_x$, and this set induces a
$K_4$ in $G[B_x]$.
These sets cover all edges of $H$, $H^*$, and
$\operatorname{Rad}(H)$, and hence all edges of $G$.
Since $G$ has no isolated vertices, every vertex belongs to
at least one bag.

For a vertex $w$ of $H$, the nodes whose skeletons contain $w$
induce a connected subtree of $\mathcal{T}^{\circ}$: vertices are identified
only across adjacent skeletons during reconstruction, so
$w$ appears in the skeleton of every node along the path
between any two such nodes. The same argument applied to
the dual SPQR tree, by Lemma~\ref{lem:spqr-duality}, proves the corresponding
assertion for a face vertex $y_f$. Thus the bags containing
any fixed vertex induce a nonempty subtree, proving \textup{(i)}.

\smallskip
\noindent
\textup{(ii)}
Reconstruct the skeletons in the two components of
$\mathcal T^\circ-xy$, retaining $e_x,e_y$, and denote the
resulting plane graphs by $H_x,H_y$, respectively.
Then
\[
V(H_x)\cap V(H_y)=\{u,v\},
\qquad
H=(H_x-e_x)\cup(H_y-e_y).
\]
Only the faces incident with $e_x$ and $e_y$ are identified
across this split, merging in pairs to give $f,g$.
Hence $B_x\cap B_y=\{u,v,y_f,y_g\}$.
Both $u$ and $v$ are incident with $f$ and $g$, giving the
four stated edges.

We claim that $H_x-e_x$ contains the edge $uv$ if and only
if $x$ is a $P$-node in $U$.
If $x$ is an $S$- or $R$-node, its skeleton is simple,
so $e_x$ is its only edge with ends $u,v$.
The subgraph attached at any other virtual edge $ab$ meets
$J_x$ only at $a,b$.
Since $\{a,b\}\ne\{u,v\}$, that subgraph cannot contain both
$u$ and $v$.
Thus $uv\notin E(H_x-e_x)$.

If $x$ is a $P$-node, each of its other core neighbours
is an $S$- or $R$-node.
The preceding argument at each such neighbour shows that
none of the attached subgraphs contains an edge $uv$.
Thus $uv$ survives precisely when $J_x$ has a real edge,
that is, when $x\in U$.
This proves the claim.
Applying it to both sides of the split shows that
$uv\in E(H)$ if and only if $x$ or $y$ is a $P$-node in $U$.
Since $G[V(H)]=H$, the same equivalence holds in $G$.
Applying this argument to $H^*$ and using
Lemma~\ref{lem:spqr-duality} proves the assertion about $y_fy_g$.

\smallskip
\noindent
\textup{(iii)}
Fix an $R$-node $\rho$.
Identify each vertex of $J_\rho$ with its corresponding
vertex of $H$, and each face $f'$ of $J_\rho$ with the
vertex $y_f\in F_\rho$ corresponding to its image face $f$.
This gives a bijection from
$V(\operatorname{Rad}(J_\rho))$ to $B_\rho$.
The boundary of $f$ is obtained from that of $f'$ by
replacing virtual edges with paths whose internal vertices
lie outside $J_\rho$.
Thus the two boundaries contain exactly the same vertices
of $J_\rho$.
Consequently, the edges of $G[B_\rho]$ between
$V(J_\rho)$ and $F_\rho$ are precisely the radial edges
under this bijection.

Consider an edge $uv$ of $G[B_\rho]$ with $u,v\in V(J_\rho)$.
Then $uv\in E(H)$.
If it is not a real edge of $J_\rho$, it belongs to the
subgraph replacing a virtual edge $ab$.
That subgraph meets $J_\rho$ only at $a,b$, so
$\{u,v\}=\{a,b\}$.
Thus $u,v$ are opposite vertices of the radial face
corresponding to an edge of $J_\rho$.
The dual argument gives the same conclusion for edges
with both ends in $F_\rho$.
All remaining edges of $G[B_\rho]$ are therefore diagonals
of radial faces.

If $e=uv$ is real, its two incident skeleton faces correspond
to its incident faces $f,g$ in $H$.
Hence both $uv$ and $y_fy_g$ belong to $G[B_\rho]$.
If $e$ is virtual and corresponds to a core edge $\rho y$,
part \textup{(ii)} shows that the primal diagonal is present
only if $y$ is a $P$-node in $U$, and the dual diagonal
only if $y$ is an $S$-node in $U$.
These conditions are mutually exclusive.
\end{proof}

\begin{corollary}\label{cor:bag-counts}
Every clique of $G$ is contained in some $G[B_x]$.
Moreover,
\begin{align}\label{eq:bag-clique-counts}
N(G,K_4)
&=\sum_{x\in V(\mathcal T^\circ)}N(G[B_x],K_4)
-|E(\mathcal T^\circ[U])|,\quad
N(G,K_5)
=\sum_{x\in V(\mathcal T^\circ)}N(G[B_x],K_5).
\end{align}
\end{corollary}

\begin{proof}
Let $C$ be a clique of $G$.
By Lemma~\ref{lem:vertex-face-bags}(i), for each $v\in V(C)$,
the nodes $x$ with $v\in B_x$ induce a nonempty
subtree of $\mathcal T^\circ$.
These subtrees pairwise intersect, since every edge
of $C$ has both ends in some bag.
By the Helly property for subtrees of a tree,
\[
\bigcap_{v\in V(C)}
\{x\in V(\mathcal T^\circ):v\in B_x\}
=
\{x\in V(\mathcal T^\circ):V(C)\subseteq B_x\}
\ne\emptyset.
\]
Thus $C$ is contained in some $G[B_x]$.
If $V(C)\subseteq B_x\cap B_y$, then
$V(C)\subseteq B_z$ for every node $z$ on
the path from $x$ to $y$ in $\mathcal T^\circ$
by Lemma~\ref{lem:vertex-face-bags}(i).
Hence the nodes whose bags contain $C$
induce a nonempty subtree of $\mathcal T^\circ$.

Let $C$ be a copy of $K_4$ contained in exactly
$t$ bags. The nodes whose bags contain $C$
induce a subtree of $\mathcal T^\circ$ with
$t$ nodes and $t-1$ edges. Thus there are exactly $t-1$ edges $xy\in E(\mathcal T^\circ)$
with $V(C)\subseteq B_x\cap B_y$.
By Lemma~\ref{lem:vertex-face-bags}(ii),
$G[B_x\cap B_y]\cong K_4$ if and only if
$xy\in E(\mathcal T^\circ[U])$.
Hence the total number of repeated counts in $\sum_{x\in V(\mathcal T^\circ)}N(G[B_x],K_4)$
is $|E(\mathcal T^\circ[U])|$, and the first
identity follows.

Every copy of $K_5$ is contained in at least one
bag, as shown above. Suppose that a copy $C$ of
$K_5$ is contained in two distinct bags.
Since the nodes whose bags contain $C$ induce
a subtree of $\mathcal T^\circ$, this subtree
contains an edge $xy$.
Thus $V(C)\subseteq B_x\cap B_y$, contradicting
$|B_x\cap B_y|=4$ by Lemma~\ref{lem:vertex-face-bags}(ii).
Hence every copy of $K_5$ is contained in
exactly one bag, proving the second identity.
\end{proof}

Recall that, for an $S$- or $P$-node $x$, the number $z_x$
counts the triangles of $H$ and $H^*$ with home $x$, and
$z_x\in\{0,1\}$ by Lemma~\ref{lem:triangle-home}.
For $x\in U\cup W$, let $a_x$ be the number of edges of
$\mathcal T^\circ[U]$ incident with $x$; in particular,
$a_x=0$ for $x\in W$.

\begin{lemma}\label{lem:SP-clique-counts}
For every $x\in U\cup W$,
\[
N(G[B_x],K_4)=q_x(1+a_x)+2z_x,
\qquad
N(G[B_x],K_5)=q_xz_x.
\]
Moreover, if $x\in U$ and $c_x=2$, then $z_x=1$ if and
only if $a_x=2$.
\end{lemma}

\begin{proof}
By duality, assume that $x$ is an $S$-node.
Write $C:=J_x$, and let $f,g$ be the faces of $H$
corresponding to the two faces of $C$.
Then $B_x=V(C)\cup\{y_f,y_g\}$.
Each subgraph substituted for a virtual edge meets $C$
only at the ends of that edge, so $G[V(C)]$ is a spanning
subgraph of $C$.
Since $C$ has at most one real edge, every vertex of $C$
is incident with a virtual edge.
Lemma~\ref{lem:vertex-face-bags}(ii) therefore implies that
both $y_f$ and $y_g$ are adjacent to every vertex of $C$.

Every real edge of $C$ belongs to $G[V(C)]$.
For a virtual edge $uv$ of $C$, its ends are adjacent in
$G$ if and only if the corresponding core neighbour is
a $P$-node in $U$.
Thus the edges of $G[V(C)]$ correspond precisely to the
$H$-active edges of $C$.
Since $x$ has no $S$-node neighbour,
Lemma~\ref{lem:vertex-face-bags}(ii) also gives
\[
y_fy_g\in E(G)
\quad\Longleftrightarrow\quad x\in U
\quad\Longleftrightarrow\quad q_x=1.
\]
When $q_x=1$, it follows that $G[V(C)]$ has $a_x+1$ edges.

By Lemma~\ref{lem:triangle-home}, $G[V(C)]$ contains exactly
$z_x$ triangles: a triangle occurs precisely when $d_x=3$
and all three skeleton edges are $H$-active.
The node corresponding to $x$ in the dual tree is a
$P$-node, so no triangle of $H^*$ has home $x$.
As a subgraph of a cycle, $G[V(C)]$ contains no $K_4$.

Every $K_4$ in the bag consequently contains one or both
of $y_f,y_g$.
There are $2z_x$ copies containing exactly one.
A copy containing both exists only when $q_x=1$ and is
then determined by an edge of $G[V(C)]$.
This gives $N(G[B_x],K_4)=q_x(a_x+1)+2z_x$.
Similarly, a $K_5$ consists of $y_f,y_g$ and a triangle
of $G[V(C)]$, giving $N(G[B_x],K_5)=q_xz_x$.

Finally, if $x\in U$ and $c_x=2$, then $d_x=3$.
The graph $G[V(C)]$ has three vertices and $a_x+1$ edges,
so it is a triangle if and only if $a_x=2$.
\end{proof}

\begin{lemma}\label{lem:SP-four-deficits}
For every $x\in U\cup W$, it holds
\[
N(G[B_x],K_4)\le4(d_x-2)+\frac{a_x}{2}.
\]
If $c_x=2$, then
\[
N(G[B_x],K_4)=
\begin{cases}
1, & a_x=0,\\
2, & a_x=1,\\
5, & a_x=2.
\end{cases}
\]
If $c_x\ge3$, then
\[
N(G[B_x],K_4)\le4(d_x-2)+\frac{a_x}{2}-2.
\]
In particular, if $c_x=2$ and $x$ has an $R$-node neighbour
in $\mathcal T^\circ$, then
$N(G[B_x],K_4)\le(a_x+3)/2$.
\end{lemma}

\begin{proof}
By Lemma~\ref{lem:SP-clique-counts},
$N(G[B_x],K_4)=q_x(a_x+1)+2z_x$.
If $q_x=0$, then $x\in W$, $a_x=0$, and $d_x=c_x\ge3$.
Hence
\[
N(G[B_x],K_4)=2z_x\le2\le4(d_x-2)-2.
\]

Now suppose that $q_x=1$, so $x\in U$ and
$0\le a_x\le c_x=d_x-1$.
If $c_x=2$, then $d_x=3$, and the three exact values follow
from Lemma~\ref{lem:SP-clique-counts}.
Each satisfies the first bound.
If $c_x\ge3$, then $d_x\ge4$, so
\[
N(G[B_x],K_4)
=a_x+1+2z_x
\le d_x+2
\le4(d_x-2)+\frac{a_x}{2}-2.
\]

Finally, if $c_x=2$ and $x$ has an $R$-node neighbour,
then $x\in U$ and $a_x\le1$.
Thus $z_x=0$ and
\[
N(G[B_x],K_4)=a_x+1\le\frac{a_x+3}{2}.\qedhere
\]
\end{proof}

\begin{lemma}\label{lem:R-bag-deficits}
Suppose that $\mathcal T^\circ$ has at least two nodes,
and let $\rho\in\mathcal R$.
Then $G[B_\rho]$ is $K_5$-free.
Moreover:
\begin{enumerate}
\item[(i)]
If $J_\rho\not\cong K_4$, then $d_\rho\ge8$ and
\(
N(G[B_\rho],K_4)\le3d_\rho-c_\rho-8.
\)
\item[(ii)]
If $J_\rho\cong K_4$, then
\[
N(G[B_\rho],K_4)\le
\begin{cases}
12, & c_\rho=1,\\
12-c_\rho, & c_\rho\ge2.
\end{cases}
\]
If, in addition, $\rho$ has an $R$-node neighbour in
$\mathcal T^\circ$, then $N(G[B_\rho],K_4)\le9$.
\end{enumerate}
\end{lemma}
\begin{proof}
By Lemma~\ref{lem:vertex-face-bags}, $G[B_\rho]$ is a spanning
subgraph of $\widehat{\Rad(J_\rho)}$.
By Lemma~\ref{lem:R-full-cliques},
$\widehat{\Rad(J_\rho)}$ is $K_5$-free,
so $G[B_\rho]$ is $K_5$ free.

We next prove \textup{(i)}. We first count the copies of $K_4$ in $G[B_\rho]$
with two vertices of $J_\rho$ and two face vertices. Let $\{u,v,f',g'\}$ be the vertex set of a $K_4$ in
$\widehat{\Rad(J_\rho)}$, where $u,v\in V(J_\rho)$ and
$f',g'$ are faces of $J_\rho$.
By construction, $uv\in E(J_\rho)$ if and only if
$uv\in E(\widehat{\Rad(J_\rho)})$.
Thus $uv\in E(J_\rho)$, and $u,v$ lie on the boundaries of both $f'$ and $g'$.
Since $J_\rho$ is 3-connected, $J_\rho-\{u,v\}$ is connected.
By Lemma~\ref{lem:common-faces}, $f',g'$ are the two faces incident with $uv$.
Let $f,g\in F(H)$ denote the faces corresponding to $f',g'$,
respectively.
If $uv$ is real, then $uv$ lies on the boundaries of $f,g$ in $H$.
Hence $uv\in E(H)$ and $fg\in E(H^*)$,
so $uv,y_fy_g\in E(G[B_\rho])$.
If $uv$ is virtual, Lemma~\ref{lem:vertex-face-bags}(ii)
implies that at most one of $uv,y_fy_g$ belongs to $E(G[B_\rho])$.
Thus $G[B_\rho]$ contains exactly $d_\rho-c_\rho$ copies
of $K_4$ with two vertices of $J_\rho$ and two face vertices.
Consequently, $N(G[B_\rho],K_4)
\le N(\widehat{\Rad(J_\rho)},K_4)-c_\rho.$

If $J_\rho\not\cong K_4$, then $J_\rho$ has at least five
vertices, all of degree at least three.
Hence $2d_\rho\ge15$, so $d_\rho\ge8$.
By Lemma~\ref{lem:R-full-cliques}, $N(G[B_\rho],K_4)\le3d_\rho-c_\rho-8.$

We now prove \textup{(ii)}. Now suppose that $J_\rho\cong K_4$.
Each triangle in $G[V(J_\rho)]$ bounds a unique face of $J_\rho$.
The triangle and the corresponding face vertex give a $K_4$
in $G[B_\rho]$.
Thus the number of copies of $K_4$ in $G[B_\rho]$
containing at least three vertices of $J_\rho$ is $N(G[V(J_\rho)],K_3)+N(G[V(J_\rho)],K_4).$
Since $G[V(J_\rho)]$ has four vertices, this sum is at most
$5$, $2$, or $1$ when $G[V(J_\rho)]$ has at most
$6$, $5$, or $4$ edges, respectively.
Suppose that $G[F(J_\rho)]$ has at most
$6$, $5$, or $4$ edges.
By duality, $G[B_\rho]$ contains at most $5$, $2$, or $1$
copies of $K_4$ with at least three vertices in $F(J_\rho)$,
respectively.

If $c_\rho=1$, then $J_\rho$ has one virtual edge $uv$
and five real edges.
Let $f,g\in F(H)$ denote the faces corresponding to
the two faces incident with $uv$ in $J_\rho$.
By Lemma~\ref{lem:vertex-face-bags}(ii), at most one of
$uv,y_fy_g$ belongs to $E(G[B_\rho])$.
Since $G[V(J_\rho)]$ and $G[F(J_\rho)]$
each have four vertices, at least one of these two graphs
has at most five edges. Moreover, $G[B_\rho]$ contains exactly $d_\rho-c_\rho=5$
copies of $K_4$ with two vertices in $V(J_\rho)$
and two in $F(J_\rho)$.
Hence $N(G[B_\rho],K_4)\le5+2+5=12.$

Suppose next that $c_\rho\ge2$.
Since $J_\rho\cong K_4$, $J_\rho$ has $6-c_\rho$ real edges
and $c_\rho$ virtual edges.
For each edge $uv$ of $J_\rho$, let $f,g\in F(H)$
denote the faces corresponding to the two faces of $J_\rho$
incident with $uv$.
If $uv$ is real, both $uv$ and $y_fy_g$ belong to $E(G[B_\rho])$.
If $uv$ is virtual, at most one of $uv,y_fy_g$ belongs to
$E(G[B_\rho])$ by Lemma~\ref{lem:vertex-face-bags}(ii).
Thus $|E(G[V(J_\rho)])|+|E(G[F(J_\rho)])|
\le 2(6-c_\rho)+c_\rho=12-c_\rho\le10.$ If one of the two graphs has at most four edges,
there are at most $1+5=6$ copies of $K_4$ in $G[B_\rho]$
with at least three vertices in $V(J_\rho)$
or at least three in $F(J_\rho)$.
Otherwise, both graphs have five edges,
and there are at most $2+2=4$ such copies.
Hence $N(G[B_\rho],K_4)\le(6-c_\rho)+6=12-c_\rho.$

Finally, suppose that $\rho$ has an R-node neighbour.
Let $uv$ be the virtual edge corresponding to this neighbour,
and let $f,g\in F(H)$ denote the faces corresponding to
the two faces of $J_\rho$ incident with $uv$.
By Lemma~\ref{lem:vertex-face-bags}(ii),
$uv,y_fy_g\notin E(G[B_\rho])$.
Thus $G[V(J_\rho)]$ and $G[F(J_\rho)]$
each have at most five edges. Hence $G[B_\rho]$ contains at most four copies of $K_4$
with at least three vertices in either $V(J_\rho)$ or $F(J_\rho)$.
Since $c_\rho\ge1$, $N(G[B_\rho],K_4)\le(6-c_\rho)+4\le9.$
\end{proof}

If $\mathcal T^\circ$ has only one node $\rho$,
then by Lemma~\ref{lem:core-leaves}, $\rho$ is an R-node.
Thus $H=J_\rho$ and $G=\widehat{\Rad(J_\rho)}$.
By Lemma~\ref{lem:R-full-cliques}, $N(G,K_5)=0$.
Since $U=\emptyset$, $N(G,K_5)=\sum_{x\in U}z_x=0$. Now suppose that $\mathcal T^\circ$ has at least two nodes.
By Lemma~\ref{lem:R-bag-deficits},
$N(G[B_\rho],K_5)=0$ for every $\rho\in\mathcal R$.
By Lemma~\ref{lem:SP-clique-counts},
$N(G[B_x],K_5)=z_x$ for $x\in U$,
and $N(G[B_x],K_5)=0$ for $x\in W$.
By Corollary~\ref{cor:bag-counts},
\begin{equation}
\label{eq:five-home-count}
N(G,K_5)=\sum_{x\in U}z_x.
\end{equation}

We now determine the maximum numbers of copies of $K_4$
and $K_5$ and characterize the extremal graphs.

\begin{proposition}\label{prop:small-orders}
For $10\le n\le13$,
\[
N(G,K_4)\le
\begin{cases}
16, & n=10,\\
19, & n=11,\\
20, & n=12,\\
25, & n=13,
\end{cases}
\]
with equality if and only if $G\cong\widehat Q$, where $Q$
is $\operatorname{Rad}(W_4)$,
$\operatorname{Rad}(C_3\square K_2)$,
$\operatorname{Rad}(W_5)$, or a graph $G_1$ from
Construction~$\ref{con}$, respectively.
\end{proposition}

\begin{proof}
Recall that $\mathcal R$ is the set of R-nodes
in $\mathcal T^\circ$.
By Lemma~\ref{lem:core-leaves} and \eqref{eq:spqr-size},
$1\le|\mathcal R|\le2$, since $n\le13$.

Suppose first that $|\mathcal R|=2$.
Then $\mathcal T^\circ$ is a path whose ends are the two R-nodes.
Since every R-skeleton other than $K_4$ has at least eight edges,
\eqref{eq:spqr-size} implies that both R-skeletons are $K_4$
and the path has $n-12$ internal nodes.
Thus $n=12$ or $13$. If $n=12$, the two R-nodes are adjacent. Since $J_\rho\cong K_4$, we have
$\widehat{\operatorname{Rad}(J_\rho)}\cong K_{2,2,2,2}$,
which contains $16$ copies of $K_4$.
By Lemma~\ref{lem:vertex-face-bags},
exactly one face of $\operatorname{Rad}(J_\rho)$
has no diagonals added.
In $K_{2,2,2,2}$, each of the two diagonals of this face
belongs to four copies of $K_4$, and exactly one copy
contains both diagonals.
Thus $N(G[B_\rho],K_4)=16-4-4+1=9.$
By Corollary~\ref{cor:bag-counts}, $N(G,K_4)=2\cdot9=18.$ If $n=13$, the path has one internal node. For each $\rho\in\mathcal R$, Lemma~\ref{lem:vertex-face-bags}
shows that $G[B_\rho]$ is obtained from $K_{2,2,2,2}$
by deleting one edge. Hence $N(G[B_\rho],K_4)=16-4=12$. The internal node $x$ has $a_x=0$.
By Lemma~\ref{lem:SP-clique-counts},
$N(G[B_x],K_4)=1$.
Hence $N(G,K_4)=2\cdot12+1=25.$ By Proposition~\ref{prop:radial-characterization},
$\operatorname{Rad}(H)\cong Q$ for some $Q\in\mathcal G^*$.
Hence $G\cong\widehat Q$.

Now suppose that $|\mathcal R|=1$.
Then $H$ is a simple $3$-connected plane graph
with $n-2$ edges. {Since $H$ and $H^*$ are simple and have at least eight edges, both have at least five vertices.}
If $n=10$, then $H$ has five vertices and eight edges.
Its degree sequence is $(4,3,3,3,3)$,
so $H\cong W_4$.
By \eqref{eq:R-full-four}, $N(G,K_4)=16.$ If $n=11$, then $|E(H)|=9$.
Since $H$ and $H^*$ are simple,
By Euler's formula,
$\{|V(H)|,|V(H^*)|\}=\{5,6\}$. By plane duality,
$\operatorname{Rad}(H)\cong\operatorname{Rad}(H^*)$.
Thus we may assume that $|V(H)|=6$ and $|V(H^*)|=5$.
Hence $H^*\cong K_5-e$ and
$H\cong C_3\square K_2$.
By \eqref{eq:R-full-four}, $N(G,K_4)=19.$

If $n=12$, then $|E(H)|=10$.
Since a simple plane graph with five vertices has at most
nine edges, both $H$ and $H^*$ have at least six vertices.
By Euler's formula, $|V(H)|=|V(H^*)|=6$. The degree sequence of $H$ is either
$(5,3,3,3,3,3)$ or $(4,4,3,3,3,3)$. If $H$ has a vertex of degree five, then $H\cong W_5$.
Since $W_5\cong W_5^*$, by \eqref{eq:R-full-four}, $N(G,K_4)=20.$ Suppose now that $H$ has no vertex of degree five.
Then both $H$ and $H^*$ have degree sequence
$(4,4,3,3,3,3)$ and contain no copy of $K_4$. By duality, $H$ and $H^*$ each have four triangular faces.
By \eqref{eq:R-full-four}, $N(G,K_4)=10+4+4=18.$

Finally, if $n=13$, by Lemma~\ref{lem:R-full-cliques}, $N(G,K_4)\le 3(n-2)-10=23<25.$
\end{proof}

\begin{proposition}\label{prop:four-bound}
For $n\ge14$, $N(G,K_4)\le4n-29$, with equality if and
only if $G\cong\widehat Q$ for some $Q\in\mathcal G^*$
with $n$ vertices.
\end{proposition}
\begin{proof}
If $\mathcal T^\circ$ has only one node $\rho$, then
$\rho$ is an R-node and $G=\widehat{\Rad(J_\rho)}$.
By \eqref{eq:spqr-size}, $d_\rho=n-2\ge12$.
By Lemma~\ref{lem:R-full-cliques}, $N(G,K_4)\le3d_\rho-10=3n-16<4n-29.$
Hence we may assume that $\mathcal T^\circ$ has at least
two nodes, so $|\mathcal R|\ge2$. Since $\sum\limits_{x\in U\cup W}a_x=2|E(\mathcal T^\circ[U])|$,
by Corollary~\ref{cor:bag-counts},
\[
N(G,K_4)
=\sum_{\rho\in\mathcal R}N(G[B_\rho],K_4)
+\sum_{x\in U\cup W}
\left(N(G[B_x],K_4)-\frac{a_x}{2}\right).
\]
By Lemma~\ref{lem:R-bag-deficits},
if $J_\rho\cong K_4$, then $d_\rho=6$ and $N(G[B_\rho],K_4)\le12=4(d_\rho-2)-4.$
If $J_\rho\not\cong K_4$, then $d_\rho\ge8$ and $c_\rho\ge1$, so $N(G[B_\rho],K_4)
\le3d_\rho-c_\rho-8
<4(d_\rho-2)-4.$ For every $x\in U\cup W$, by Lemma~\ref{lem:SP-four-deficits},
$N(G[B_x],K_4)-\frac{a_x}{2}\le4(d_x-2).$ Thus $N(G,K_4)
\le4\sum_{x\in V(\mathcal T^\circ)}(d_x-2)-4|\mathcal R|.$
By \eqref{eq:spqr-size}, $N(G,K_4)\le4n-16-4|\mathcal R|.$

Suppose first that $|\mathcal R|\ge3$.
If some $\rho\in\mathcal R$ satisfies $c_\rho\ge2$,
then, by Lemma~\ref{lem:R-bag-deficits}, $N(G[B_\rho],K_4)\le4(d_\rho-2)-6.$
Hence $N(G,K_4)
\le4n-16-6-4(|\mathcal R|-1)
\le4n-30.$
Otherwise, every R-node is a leaf of $\mathcal T^\circ$.
The tree therefore has at least three leaves,
so some S- or P-node $x$ satisfies $c_x\ge3$.
By Lemma~\ref{lem:SP-four-deficits}, $N(G[B_x],K_4)-\frac{a_x}{2}
\le4(d_x-2)-2.$
Thus $N(G,K_4)\le4n-16-4|\mathcal R|-2\le4n-30.$

Now suppose that $|\mathcal R|=2$.
By Lemma~\ref{lem:core-leaves},
$\mathcal T^\circ$ is a path whose leaves are the two R-nodes.
Every internal node $x$ satisfies $c_x=2$.
Since $d_x\ge3$ and $q_x\le1$,
$d_x=3$ and $q_x=1$, so $x\in U$. If the path has no internal node,
the two R-nodes are adjacent. For every $\rho\in\mathcal R$, by Lemma~\ref{lem:R-bag-deficits}, $N(G[B_\rho],K_4)\le4(d_\rho-2)-7.$ By \eqref{eq:spqr-size}, $\sum_ {\rho\in\mathcal R}(d_\rho-2)=n-4.$
Hence
\[
N(G,K_4)
=\sum_{\rho\in\mathcal R}N(G[B_\rho],K_4)
\le4\sum_{\rho\in\mathcal R}(d_\rho-2)-14
=4n-30.
\]
If the path has only one internal node $x$,
both neighbours of $x$ are R-nodes, so $a_x=0$.
By Lemma~\ref{lem:SP-clique-counts}, $N(G[B_x],K_4)=1=4(d_x-2)-3$.
If both R-skeletons were isomorphic to $K_4$, by \eqref{eq:spqr-size}, $n=13$.
Since $n\ge14$, some $\rho\in\mathcal R$
satisfies $J_\rho\not\cong K_4$.
By Lemma~\ref{lem:R-bag-deficits}, $N(G[B_\rho],K_4)
\le3d_\rho-c_\rho-8
\le4(d_\rho-2)-9,$
where $d_\rho\ge8$ and $c_\rho=1$.
Therefore $N(G,K_4)\le4n-16-9-4-3=4n-32.$
Suppose that the path has at least two internal nodes.
Each internal node $x$ adjacent to an R-node satisfies
$a_x=1$ and $d_x=3$.
By Lemma~\ref{lem:SP-clique-counts}, $N(G[B_x],K_4)-\frac{a_x}{2}
=2-\frac12
=4(d_x-2)-\frac52.$
Hence $N(G,K_4)
\le4n-16-4-4-\frac52-\frac52
=4n-29.$

Suppose that $N(G,K_4)=4n-29$.
Then $|\mathcal R|=2$, and the path has at least
two internal nodes. Every R-skeleton is isomorphic to $K_4$.
Every internal node belongs to $U$ and has one real edge
and two virtual edges in its skeleton.
Since no two S-nodes and no two P-nodes are adjacent,
the internal nodes alternate between S- and P-nodes.
By \eqref{eq:spqr-size},
the path has $n-12$ internal nodes. By Proposition~\ref{prop:radial-characterization},
$\operatorname{Rad}(H)\cong Q$ for some $Q\in\mathcal G^*$.
Hence $G\cong\widehat Q$.

Conversely, suppose that $G\cong\widehat Q$ for some
$Q\in\mathcal G^*$ with $n\ge14$ vertices.
By Proposition~\ref{prop:radial-characterization}, we may choose $H$
so that its SPQR tree has the core and skeletons described above. Each end R-node $\rho$ has one neighbour in $U$.
Since $J_\rho\cong K_4$, we have
$\widehat{\operatorname{Rad}(J_\rho)}\cong K_{2,2,2,2}$,
which contains $16$ copies of $K_4$.
By Lemma~\ref{lem:vertex-face-bags},
exactly one face of $\operatorname{Rad}(J_\rho)$
has only one diagonal added.
Since each edge of $K_{2,2,2,2}$ belongs to four copies of $K_4$, $N(G[B_\rho],K_4)=16-4=12.$ The two internal nodes adjacent to the R-nodes
have $a_x=1$.
The remaining $n-14$ internal nodes have $a_x=2$. By Lemma~\ref{lem:SP-clique-counts},
$N(G[B_x],K_4)=2$ for each internal node $x$
adjacent to an R-node, and $N(G[B_x],K_4)=5$
for each other internal node.
By Corollary~\ref{cor:bag-counts}, $N(G,K_4)=2\cdot12+2\cdot2+5(n-14)-(n-13)=4n-29.$
\end{proof}

\begin{proposition}\label{prop:five-bound}
For $n\ge10$, it holds 
\(
N(G,K_5)\le\max\{0,n-14\}.
\)
For $n\ge15$, with equality if and only if
$G\cong\widehat Q$ for some $Q\in\mathcal G^*$ with
$n$ vertices.
\end{proposition}
\begin{proof}
By \eqref{eq:five-home-count}, $N(G,K_5)=\sum_{x\in U}z_x.$
If $N(G,K_5)=0$, the upper bound holds.
Suppose now that $N(G,K_5)>0$. There exists $x\in U$
with $z_x=1$. Then $d_x=3$
by Lemma~\ref{lem:triangle-home}.
Thus $c_x=d_x-q_x=2$, and by
Lemma~\ref{lem:SP-clique-counts}, $a_x=2$. The component of $\mathcal T^\circ[U]$ containing $x$
has at least three nodes and hence at least two leaves.
Each leaf $y$ satisfies $a_y=1$, so $z_y=0$.
By Lemma~\ref{lem:core-leaves}, $\mathcal T^\circ$
has at least two leaves in $\mathcal R$, so $|\mathcal R|\ge2$.
Therefore, by \eqref{eq:spqr-size},
\[
N(G,K_5)
\le |U|-2
\le n-6-4|\mathcal R|-|W|
\le n-14.
\]
In particular, $N(G,K_5)>0$ implies $n\ge15$. This proves the upper bound.

Now suppose that $n\ge15$ and $N(G,K_5)=n-14$.
Thus $|\mathcal R|=2$, $W=\emptyset$, and $|U|=n-12$.
By \eqref{eq:spqr-size},
$d_\rho=6$ for every $\rho\in\mathcal R$ and
$d_x=3$ for every $x\in U$.
Thus both R-skeletons are isomorphic to $K_4$. By Lemma~\ref{lem:core-leaves},
every leaf of $\mathcal T^\circ$ is an R-node.
Since $|\mathcal R|=2$, $\mathcal T^\circ$ is a path
whose leaves are the two R-nodes. The path has $n-12$ internal nodes, all in $U$.
Each has one real edge and two virtual edges in its skeleton.
Since no two S-nodes and no two P-nodes are adjacent,
the internal nodes alternate between S- and P-nodes. By Proposition~\ref{prop:radial-characterization},
$\operatorname{Rad}(H)\cong Q$ for some $Q\in\mathcal G^*$.
Hence $G\cong\widehat Q$.

Conversely, suppose that $G\cong\widehat Q$ for some
$Q\in\mathcal G^*$ with $n\ge15$ vertices.
By Proposition~\ref{prop:radial-characterization}, we may choose $H$
so that its SPQR tree has the core and skeletons described above.
The two internal nodes adjacent to the R-nodes have $a_x=1$,
and every internal node not adjacent to an R-node has $a_x=2$.
By Lemma~\ref{lem:SP-clique-counts} and
\eqref{eq:five-home-count}, $N(G,K_5)=\sum_{x\in U}z_x=n-14.$
\end{proof}

\begin{proof}[\bf Proof of Theorem~\ref{thm:cliques}]
The assertions concerning $K_3$ were proved above (see Proposition~\ref{prop:small-triangles} and Theorem~\ref{thm}).
Propositions~\ref{prop:small-orders}, \ref{prop:four-bound},
and~\ref{prop:five-bound} establish the assertions concerning
$K_4$ and $K_5$.
\end{proof}

\section{Disproof of Conjecture~\ref{conj2}}
\label{sec:2planar}

We construct an infinite family of counterexamples to
Conjecture~\ref{conj2}.

\begin{proposition}\label{prop:2planar-counterexamples}
For every integer $r\ge1$, there exists a simple $2$-planar
graph $F_r$ with $5r+3$ vertices and $44r+1$ triangles.
\end{proposition}

\begin{proof}
Let $F$ be obtained from $K_7$ on the vertex set
$\{A,B,C,a,b,c,d\}$ by adding a vertex $v$ adjacent to
$A,B,C,b,c$.
Thus $F=K_8-\{va,vd\}$.
A $2$-planar drawing of this graph appears in
Schneck~\cite[Section~5.7.1, Fig.~5.8(b)]{ref45}.
We use the labelled drawing in
Figure~\ref{fig:2planar-counterexample}.
The triangles $abc$ and $vAC$ are vertex-disjoint and
uncrossed, and bound a bounded face and the unbounded
face of the drawing, respectively.   
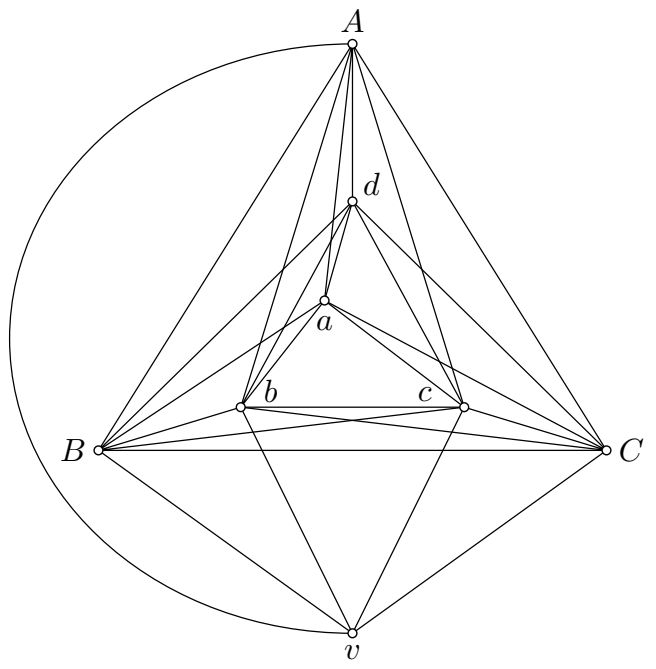
\begin{figure}[htbp]
  \centering
  \scalebox{0.48}{%
  \begin{tikzpicture}[
    x=0.70cm, y=0.70cm,
    line cap=round, line join=round,
    edge/.style={draw=black, line width=0.50pt},
    vertex/.style={circle, draw=black, fill=white, line width=0.55pt,
                   inner sep=0pt, minimum size=3.45pt},
    every node/.style={inner sep=0pt}
  ]
    \coordinate (A) at ( 0.00, 8.00);
    \coordinate (B) at (-5.00, 0.00);
    \coordinate (C) at ( 5.00, 0.00);
    \coordinate (a) at (-0.55, 2.95);
    \coordinate (b) at (-2.20, 0.85);
    \coordinate (c) at ( 2.20, 0.85);
    \coordinate (d) at ( 0.00, 4.90);
    \coordinate (v) at ( 0.00,-3.60);

    \draw[edge] (A)--(B) (A)--(C) (A)--(a) (A)--(b) (A)--(c) (A)--(d);
    \draw[edge] (B)--(C) (B)--(a) (B)--(b) (B)--(c) (B)--(d);
    \draw[edge] (C)--(a) (C)--(b) (C)--(c) (C)--(d);
    \draw[edge] (a)--(b) (a)--(c) (a)--(d) (b)--(c) (b)--(d) (c)--(d);

    \draw[edge] (v)--(B) (v)--(C) (v)--(b) (v)--(c);
    \draw[edge] (v) arc[start angle=-90, end angle=-270,
                         x radius=6.75, y radius=5.80];

    \foreach \p in {A,B,C,a,b,c,d,v}
      \node[vertex] at (\p) {};

    \tikzset{every node/.append style={scale=1.15}}
    \node[above=4pt] at (A) {$A$};
    \node[left=4pt]  at (B) {$B$};
    \node[right=4pt] at (C) {$C$};
    \node[below=5pt] at (a) {$a$};
    \node[xshift=10.2pt, yshift=5.5pt] at (b) {$b$};
    \node[xshift=-13.4pt,yshift=4.3pt] at (c) {$c$};
    \node[xshift=6.5pt,  yshift=6.0pt] at (d) {$d$};
    \node[below=4pt] at (v) {$v$};
  \end{tikzpicture}%
  }
  \caption{A $2$-planar drawing of $F=K_8-\{va,vd\}$.
    The triangles $abc$ and $vAC$ are uncrossed face boundaries.}
  \label{fig:2planar-counterexample}
\end{figure}

Since $F-v\cong K_7$ and the five neighbours of $v$
are pairwise adjacent,
\[
N(F,K_3)=\binom{7}{3}+\binom{5}{2}=45.
\]
Set $F_1:=F$ and $T_1:=abc$.
Given a drawing of $F_r$ in which $T_r$ bounds an
uncrossed bounded face, place a new copy of $F$ inside
that face and identify its outer triangle $vAC$ with
$T_r$, including the corresponding boundary edges.
Let $T_{r+1}$ be the triangle $abc$ in the new copy.

The two copies meet precisely in the identified triangle,
so the resulting graph $F_{r+1}$ is simple.
Both copies have uncrossed attachment boundaries, so
the gluing creates no new crossings.
Thus $F_{r+1}$ is $2$-planar, and $T_{r+1}$ again bounds
an uncrossed bounded face.

Each step adds five vertices.
No edge joins a new vertex to
$V(F_r)\setminus V(T_r)$, so every triangle of
$F_{r+1}$ lies entirely in $F_r$ or in the new copy.
The only triangle common to both is $T_r$.
Consequently,
\[
|V(F_{r+1})|=|V(F_r)|+5,
\qquad
N(F_{r+1},K_3)=N(F_r,K_3)+44.
\]
Starting from $|V(F_1)|=8$ and $N(F_1,K_3)=45$,
induction gives
\[
|V(F_r)|=5r+3,\qquad N(F_r,K_3)=44r+1.
\]
Finally,
\[
N(F_r,K_3)-\frac{17|V(F_r)|-49}{2}
=\frac{3r}{2}>0,
\]
so every $F_r$ violates the conjectured bound.
\end{proof}

The construction also gives
\[
\lim_{r\to\infty}\frac{N(F_r,K_3)}{|V(F_r)|}
=\frac{44}{5}.
\]
Hence any upper bound of the form
$N(G,K_3)\le c|V(G)|+O(1)$ valid for all simple
$2$-planar graphs must satisfy $c\ge44/5$.
Determining the optimal coefficient and the corresponding
extremal graphs is a natural further question.
\section*{Declaration on the use of AI}
While preparing this manuscript, generative AI was used to assist with English revision, proof presentation, \LaTeX{} formatting, and reference checking. The authors have reviewed, edited, and independently validated all AI-processed material, including mathematical statements and citations, and accept full responsibility for the final manuscript.

\end{document}